\documentclass[reqno, english]{amsart}
\usepackage{amsfonts}
\usepackage{amsthm}
\usepackage{amsmath}
\usepackage{amscd}
\usepackage[latin2]{inputenc}
\usepackage{t1enc}
\usepackage[mathscr]{eucal}
\usepackage{indentfirst}
\usepackage{graphicx}
\usepackage{graphics}
\usepackage{pict2e}
\usepackage{epic}
\usepackage{bbm}
\usepackage{bm}
\usepackage[left=1in,right=1in,top=1in,bottom=1in]{geometry}
\usepackage{epstopdf} 
\usepackage{enumitem}

\usepackage{amssymb}

\usepackage{hyperref}
\hypersetup{
	colorlinks = true,
    linkcolor = blue,
    citecolor = blue,
    urlcolor = black,
}

\theoremstyle{plain}
\newtheorem{thm}{Theorem}[section]

\newtheorem{lemma}[thm]{Lemma}
\newtheorem{cor}[thm]{Corollary}
\newtheorem{prop}[thm]{Proposition}
\newtheorem{claim}[thm]{Claim}
\newtheorem{question}[thm]{Question}
\newtheorem{Ex}[thm]{Example}

\newtheorem{obs}[thm]{Observation}
\newtheorem{conj}[thm]{Conjecture}

\newtheorem*{claim*}{Claim}

\newcommand{\K}{\mathcal K}

\newcommand{\mS}[0]{\mathfrak{S}}

\newcommand{\prr}{p}
\newcommand{\pr}{\mathbb{P}}

\newcommand{\E}[0]{\mathbb{E}}
\newcommand{\mbone}[0]{\mathbbm{1}}

\newcommand{\beq}[1]{\begin{equation}\label{#1}}
\newcommand{\enq}[0]{\end{equation}}

\newcommand{\bn}[0]{\bigskip\noindent}
\newcommand{\mn}[0]{\medskip\noindent}
\newcommand{\nin}[0]{\noindent}

\newcommand{\sub}[0]{\subseteq}
\newcommand{\sm}[0]{\setminus}

\renewcommand{\dots}[0]{,\ldots,}

\newcommand{\more}[0]{~~\mbox{\raisebox{-.9ex}{$\stackrel{\textstyle{>}}{\sim}$}}~~}

\newcommand{\A}[0]{{\mathcal A}}

\newcommand{\cee}[0]{{\mathcal C}}

\newcommand{\h}[0]{{\mathcal H}}
\newcommand{\I}[0]{{\mathcal I}}

\newcommand{\Oh}[0]{{\mathcal O}}

\newcommand{\V}[0]{{\mathcal V}}

\newcommand{\ra}[0]{\rightarrow}
\newcommand{\Ra}[0]{\Rightarrow}

\newcommand{\Lra}[0]{\Leftrightarrow}

\newcommand{\Rr}[0]{\mbox{${\bf R}$}}

\newcommand{\0}[0]{\emptyset}

\newcommand{\ga}[0]{\alpha }
\newcommand{\gb}[0]{\beta }
\newcommand{\gc}[0]{\gamma }
\newcommand{\gd}[0]{\delta }

\newcommand{\gl}[0]{\lambda }

\newcommand{\go}[0]{\omega}
\newcommand{\gO}[0]{\Omega}

\newcommand{\gs}[0]{\sigma}

\newcommand{\gz}[0]{\zeta}
\newcommand{\eps}[0]{\varepsilon }
\newcommand{\vt}[0]{\vartheta}
\newcommand{\vs}[0]{\varsigma}

\newcommand{\vp}[0]{\varphi}

\newcommand{\supp}{\text{supp}\,}

\newcommand{\ub}[0]{b}%{\underline{b}}
\newcommand{\uc}[0]{c}%{\underline{c}}
\newcommand{\ud}[0]{d}%{\underline{d}}
\newcommand{\uj}[0]{\underline{j}}
\newcommand{\uv}[0]{v}%{\underline{v}}
\newcommand{\uw}[0]{\underline{w}}

\newcommand{\uzero}[0]{\underline{0}}

\newcommand{\prh}[1][]{\pr_h}

\newcommand{\GP}[0]{G(P)}

\usepackage{xcolor}

\newcommand{\tdf}{\tilde{f}}

\newcommand{\Var}[0]{{\rm Var}}
\newcommand{\var}[0]{{\rm var}}
\newcommand{\win}[0]{{\rm win}}
\newcommand{\Win}[0]{{\rm Win}}
\newcommand{\gap}[0]{{\rm gap}}
\newcommand{\height}[0]{{\rm ht}}

\newcommand{\mbR}[0]{\mathbb{R}}
\newcommand{\mbZ}[0]{\mathbb{Z}}

\newcommand{\conv}[0]{{\rm conv}}

\newcommand{\llangle}[0]{\langle\hspace{-0.02in}\langle}
\newcommand{\rrangle}[0]{\rangle\hspace{-0.02in}\rangle}

\newcommand{\cov}[0]{~\cdot\hspace{-0.05in}\succ}
\newcommand{\recovby}[0]{<\hspace{-0.07in}\cdot~}

\newcommand{\hO}[0]{\hat{0}}
\newcommand{\hone}[0]{\hat{1}}
\newcommand{\gss}[0]{\vs}
\newcommand{\lxr}[0]{\langle x\rangle}
\newcommand{\llxr}[0]{\llangle x\rrangle}

\title{Balancing extensions in posets of large width.}

\author{Max Aires, Jeff Kahn}

\begin{document}

\maketitle
\begin{abstract}

We revisit classic balancing problems for linear extensions of a 
partially ordered set $P$, proving results that go 
far beyond many of the best earlier results on this topic.
%well beyond most of what was previously known.  
For example, 
with $\prr(x\prec y)$ the probability
that $x$ precedes $y$ in a uniform linear extension, 
$\gd_{xy}=\min\{\prr(x\prec y),\prr(y\prec x)\}$, and 
$\gd(P)=\max \gd_{xy}$, we show that
$\gd(P)$ tends to $1/2$ as $n:=|P|\ra \infty$ if 
$P$ has width $\gO(n)$ or $\go(\log n)$ minimal elements, 
and is at least $1/e-o(1)$ if $P$ has width $\go(\sqrt{n})$ or height $o(n)$.

Motivated by both consequences for balance problems
and intrinsic interest, we also consider several old and new 
parameters associated with $P$.
Here, in addition to balance, we study relations between the parameters
%and to the conjecture (proposed by Kahn and Saks in 1984)
%that $\gd(P)\ra 1/2$ as the width of $P$ tends to infinity,
and suggest various questions that are thought to be worthy of further investigation.

\end{abstract}

\section{Introduction: balance theorems}
\label{Intro}

Probabilities associated with uniform distribution on
the collection of linear extensions of a partially 
ordered set ({\em poset}) have been the subject of some
fascinating problems and results over the last 40 or more years,
and will again be our subject here.
(See e.g.\ \cite{Chan-Pak} for a recent survey.)

Following a standard abuse, we identify a poset $P$ with its set of elements.
A \emph{linear extension} of an $n$-element $P$
may then be viewed as either  
an order-preserving bijection from $P$ to $ [n]:=\{1\dots n\}$,
or as a linear ordering of the elements of $P$
respecting the poset relations.
We use $E(P)$ for the set of linear extensions of $P$ and $e(P)$ for $|E(P)|$,
and throughout the paper use $f$ or $\gss$
for a uniform member of $E(P)$,
according to whether we are thinking of a bijection or an ordering.
When thinking of $\gss$ we use (e.g.) $\prr(x\prec y)$---the probability
that $x$ precedes $y$ in $\gss$---in place of $\pr(f(x)<f(y))$.
Further terminology 
(e.g.\ chain, antichain, height $\height(P)$, and average height $h(x)$)
and usage notes are given at the end of this section.

Our first goal in this paper is to say something new 
about a few old ``balancing'' problems for linear extensions,
and in this section we mainly 
sketch what we can do in that direction, slightly postponing 
other aspects, in particular 
%the parameters of Section~\ref{Parameters}, introduction and study of which 
the parameters introduced in Section~\ref{Parameters} that
are our second main concern.
%(See the outline below for an overview.) 

Possibly the best known of all poset problems 
is the following ``1/3-2/3 Conjecture." 
Here and throughout we use $\gd_{xy}$ for $\min\{\prr(x\prec y),\prr(y\prec x)\}$,
$\gd(P)$ for $\max \gd_{xy}$, the max over distinct $x,y\in P$, and $\gd_x$ for $\max_{y\neq x}\gd_{xy}$.
%and $\gd(P)$ for $\max \gd_{xy}$, the max over distinct $x,y\in P$.

\begin{conj}
\label{C1/3}
If the finite poset P is not a chain then 
\beq{1/3}
\gd(P)\geq 1/3.
\enq
\end{conj}
\nin
\nin
This was first proposed by Kislitsyn \cite{Kislitsyn}
and again by Fredman (circa 1975) and Linial \cite{Linial},
all motivated by questions about sorting.
From an algorithmic standpoint what's wanted is existence of
\emph{some} universal $\gd>0$ 
such that $\gd(P)\geq \gd$ whenever $P$ is not a chain.

This existence was shown
in \cite{Kahn-Saks}, with $\gd=3/11$ and a proof
based especially on the Aleksandrov-Fenchel Inequalities from
the theory of mixed volumes (e.g. \cite{Schneider}).
Easier proofs with smaller $\gd$'s,
based on the Brunn-Minkowski Theorem
and an idea of Gr\"unbaum \cite{Grunbaum}, were given
(independently) in \cite{Karzanov-Khachiyan} 
and \cite{Kahn-Linial}, which show, respectively, $\gd(P) > 1/e^2$
and $\gd(P) > 1/(2e)$.
(The present work, too, will often lean on geometry;
see Section~\ref{Geometry}.)

The current best bound on $\gd$ 
is that of Brightwell, Felsner and Trotter \cite{BFT}, 
who combine the machinery of \cite{Kahn-Saks} with 
a special case of the \emph{Cross-Product Conjecture} to
improve the $3/11$ of \cite{Kahn-Saks} to $(5-\sqrt{5})/10$.
(The still open Cross-Product Conjecture is also from \cite{BFT};
see \cite{CPP} for status and recent progress.)
What's particularly striking about \cite{BFT} is that the bound on $\gd$ is 
best possible for the more general infinite version of
the problem (to which all the arguments just mentioned still apply).

There has also been a substantial amount of work proving Conjecture~\ref{C1/3}
for particular classes of posets, but we will not review these; see \cite{Brightwell}
for an account of what was known as of 25 years ago, and
\cite{Chan-Pak} for a more recent summary.

While Conjecture~\ref{C1/3} is tight, it seems likely 
that one can usually do better .
%except in quite special situations.
The following natural guess in this direction was proposed in \cite{Kahn-Saks} and 
is now known as the 
\emph{Kahn-Saks Conjecture}.
(As usual, the \emph{width}, $w(P)$, is the size of a largest antichain in $P$.)
Another obvious guess, again seemingly first suggested in \cite{Kahn-Saks},
is that $\gd(P)>1/3$ if $w(P)\geq 3$.

\begin{conj}
\label{CKS}
If $w(P) \ra \infty$, then
\beq{1/2}
\gd(P)\ra 1/2.
\enq
\end{conj}
\nin
Convergence here is meant uniformly in $P$; that is,
$\gd(P)>1/2-o(1)$, where $o(1)\ra 0$ as $w(P)\ra\infty$.
Alternatively, one may regard $P$ as an entry in a
\emph{sequence} of posets with widths tending to infinity,
but this formality seems unnecessary (and annoying) and we avoid it.

Progress on Conjecture~\ref{CKS} has been limited.  The best general result 
is the following theorem of Koml\'os \cite{Komlos}, which he derives from
a beautiful probabilistic statement to which we will return below
(see Theorem~\ref{TkKom}).
%which shows that \eqref{1/2} holds if
%$|\min (P)|=\gO(n)$ (equivalently, $|\min P| > \eps n$ for 
%$\eps=\eps_n$ tending slowly enough to zero).
\begin{thm}%[\cite{Komlos}]
\label{TKom}
A sufficient condition for \eqref{1/2} is that $|\min(P)| =\gO(n)$.
\end{thm}
\nin
(So the convergence in \eqref{1/2} is now as $n\ra \infty$.)
Of course an alternate form of Theorem~\ref{TKom} says
that \eqref{1/2} holds if $|\min (P)| > \eps n$ for 
$\eps=\eps_n$ tending slowly enough to zero.  
We will omit similar reformulations below.

In other work related to Conjecture~\ref{CKS},
Korshunov \cite{Korshunov} showed that \eqref{1/2} holds 
w.h.p.\footnote{\emph{With high probability,} meaning with 
probability tending to 1 as $n = |P| \ra\infty$.}
for random posets in the sense of \cite{Kleitman-Rothschild},
where it is shown that the height of a random poset is 3 w.h.p.
Very precise estimates implying \eqref{1/2} for Young diagrams of 
\emph{bounded} width are given in \cite{CPP2}, to which we refer for related work.

\nin

\mn

Our main results in the direction of Conjecture~\ref{CKS} are as follows.

\begin{thm}\label{TKS}
If $P$ satisfies any one of the following then \eqref{1/2} holds:

\nin
{\rm (a)}
$|\min(P)|\gg \log n$; 

\nin
{\rm (b)}
$\min_{x\in P}h(x) \gg 1$;

\nin
{\rm (c)}
$w(P) =\gO(n)$; 

\nin
{\rm (d)}
$P = C\sqcup A$ with $C$ a chain and $A$ an antichain of size $\go(1)$.
\end{thm}

\nin
(Below we employ the natural misuse of ``Theorem~\ref{TKS}(a)''---or just ``(a)'' 
if the theorem is clear---for the statement ``(a) implies \eqref{1/2}.'')

\nin
Of course each of (a) and (c) substantially strengthens 
Theorem~\ref{TKom}.
(As Koml\'os \cite{Komlos} observes, the $\log n$ in (a) is
best possible in a ``local'' sense; see the remarks following Theorem~\ref{T1/e}.)
And in view of \cite{Kleitman-Rothschild}, \cite{Korshunov}
is a special case of
\eqref{1/2} for \emph{height} 3, whereas Theorem~\ref{TKS}(c) implies
\eqref{1/2} for \emph{any fixed height}.
We include (d), which is unrelated to the other parts of Theorem~\ref{TKS}, 
both to stress that even in this very special 
case, Conjecture~\ref{CKS} is not easy,
and because we think the proof is interesting and possibly suggestive.
%(See also the remarks following Theorems~\ref{TKS'} and \ref{T1/e}.)

The proofs of (a)-(c)
%, combined with a generalization of that theorem (Theorem~\ref{TkKom} below)
give more.
For $X\sub P$ define $\gd_k(X)$ to be the largest $\gd$ for which there is $Y\sub X$ of size $k$
such that $\gss|_Y$ takes each of its $k!$ 
possible values with probability at least $\gd$.  
\begin{thm}\label{TKS'}
For fixed k, if P satisfies (a), (b) or (c) of Theorem~\ref{TKS}, then
\beq{1/k!}
\gd_k(P)\ra 1/k!.
\enq
\end{thm}
\nin
On the other hand, (d) of Theorem~\ref{TKS}
(\emph{a fortiori} the $w(P)\ra\infty$ of Conjecture~\ref{CKS}) doesn't even give \eqref{1/k!}
with $k=3$; e.g.\ it fails with $m=2^t$, $A=\{x_1\dots x_t\}$, $C=\{y_1<\cdots < y_m\}$, and 
the additional relations $x_i\recovby y_{2^i}$.

\nin
\emph{Remark.}  Theorems~\ref{TKS}/\ref{TKS'}(b) and Theorem~\ref{TKS}(d)
differ from our other balance results
%---and any others we know of---
in assuming
only that something is large in absolute terms, 
rather than as a function of $n$.

\mn

In much of what follows we will be interested in a conclusion between \eqref{1/3}
and \eqref{1/2}, namely
\beq{1/e}
\gd(P) > 1/e-o(1).
\enq
The $1/e$, corresponding to the fact 
%\cite{Grunbaum}
that at least $1/e$ of the volume of a convex body 
$K\sub \Rr^d$ lies on each side of any hyperplane
through its centroid 
(see Gr\"unbaum's Theorem~\ref{TGrunbaum} below), 
is natural in the present context, suggesting that one
try to prove Conjecture~\ref{CKS} with \eqref{1/2} replaced by \eqref{1/e}.
In particular this would contain the intersection of 
Conjectures~\ref{C1/3} and \ref{CKS}; that is, 1/3-2/3 for large width.

In \cite{Friedman} Friedman
used variations on the geometry underlying the $1/(2e)$ bound 
of \cite{Kahn-Linial} to give the following general sufficient conditions for \eqref{1/e}.
Here we use $\ga(x)=|\{y:y\leq x\}|$ and $\gb(x)=|\{y:y\geq x\}|$.

\begin{thm}\label{TF}
If $P$ satisfies any one of the following then \eqref{1/e} holds:

\nin
{\rm (a)}
$|\min(P)|  \gg \sqrt{n}$;

\nin
{\rm (b)}
$\height(P)<2\log_2\log n-\go(1)$.

\nin
{\rm (c)}
$\min_{f\in E(P)}\sum \left[f(x)/\ga(x)+ (n+1-f(x))/\gb(x)\right]\gg n$.

\end{thm}

\nin
\emph{Notes.}
The main condition here is (c), from which Friedman derives the other two.
(In \cite{Friedman} the min in (c) is over \emph{all} bijections $f:P\ra[n]$, but 
this is easily seen to be equivalent.)
%, i.e.\ the min is achieved with $f$ an extension.)

Theorem~\ref{TF}(a), like Theorem~\ref{TKom},
is subsumed by Theorem~\ref{TKS},
and each of (b),(c) is greatly improved by the corresponding part of the next result,
which summarizes what we can say about sufficient 
conditions for \eqref{1/e}.
Recall that a \emph{fractional matching} of $\h\sub 2^V$ is a
$\gl: \h\ra [0,1]$ satisfying
$\sum_{x\in A\in\h}\gl_A \leq 1 ~\forall x\in V$.

\begin{thm}\label{T1/e}
If $P$ satisfies any one of the following then \eqref{1/e} holds:

\nin
{\rm (a)}
$w(P)\gg \sqrt{n}$;

\nin
{\rm (b)}
$\height(P) \ll n$;

\nin
{\rm (c)}
$\sum h(x)/\ga(x) \gg n$;
%$\sum_{x\in P} h(x)/|\langle x\rangle| \gg n$;

\nin
{\rm (d)}
$\log e(P) \gg n$;

\nin
{\rm (e)}
there is a fractional matching $\gl$ of $\A=\A(P)$ satisfying
\beq{gd.frac.m}
\mbox{$\sum_{A\in \A}\gl_A |A|^2\gg n.$}
\enq
\end{thm}
\nin
Here part (c) implies not only 
Theorem~\ref{TF}(c), but that \eqref{1/e} holds if the condition in Theorem~\ref{TF}(c)
holds \emph{on average}; for in that case
one of
$\sum h(x)/\ga(x)$, 
$\sum (n+1-h(x))/\gb(x)$ is $\go(n)$,
and \eqref{1/e} follows from Theorem~\ref{T1/e}(c)
applied to either $P$ or its dual, $P^*$ (since, trivially, $\gd_{P^*} =\gd_P$).

\nin

\nin
\emph{Remarks.}
As we will see, (a)-(c) of Theorems~\ref{TKS}/\ref{TKS'},
and
(a) and (e) of Theorem~\ref{T1/e}, are \emph{local} in the sense that
the balanced pairs (or $k$-tuples in Theorem~\ref{TKS'})
can be chosen from $\min(P)$ in Theorems~\ref{TKS}/\ref{TKS'}(a,b);
$\cup\{A:A\in\supp(\gl)\}$ in 
Theorem~\ref{T1/e}(e); and any antichain of size $\gO(w(P))$ in the other two cases.

At this local level, the bounds of Theorem~\ref{TKS}(a)
(so also Theorem~\ref{TKS'}(a))
and  Theorem~\ref{T1/e}(a) are essentially best possible;
for Theorem~\ref{TKS}(a) this is shown
by Koml\'os' example \cite{Komlos} of $P$ consisting of $\Theta(\log n)$ unrelated chains of 
geometrically increasing lengths, or by the example in the paragraph
following Theorem~\ref{TKS'};
for Theorem~\ref{T1/e}(a) see Example~\ref{MaxT4.1} below.

In contrast, Theorem~\ref{TKS}(d) is \emph{non}local
(and does not have a Theorem~\ref{TKS'} version);
e.g.\ when $C=\{y_1<\cdots <y_m\}$, with $m=2^k$,
and $A=\{x_1\dots x_k\}$ with $x_i\recovby y_{2^i}$,
where a balanced pair must be of the form $\{x_i,y_j\}$.

\bn
\textbf{Outline.}

Section~\ref{Parameters}
introduces several new parameters 
associated with a poset.
These are,
as mentioned above, the paper's second main concern, and
Section~\ref{Parameters} serves as a sort of second introduction, outlining 
what we will say about the parameters themselves and their relevance to balance questions.
Highlights here (a couple of whose proofs are deferred to \cite{AK2})
include:
Theorems~\ref{Twwin}-\ref{Ttaugs}
and the stronger versions of the first two in Theorem~\ref{Tsumwin} and
Conjecture~\ref{A1};
Theorem~\ref{Twsig} and Corollary~\ref{Cpi}, which
give conditions implying \eqref{1/2} for $P$ of \emph{bounded} width and, with
Theorem~\ref{Ttaugs}, imply the surprising equivalence of Conjecture~\ref{CKS}
and the seemingly stronger Conjecture~\ref{maxC}; 
and Theorem~\ref{Twin}, which 
connects \eqref{1/e} to the parameter ``$\win$'' treated in Theorems~\ref{Twwin} and \ref{Tsumwin}.

Sections~\ref{Geometry} and \ref{Log-concavity} provide the connection to geometry and 
briefly recall log-concavity background,
and 
Section~\ref{Consequences} derives a few assertions---especially
Theorem~\ref{Tsumwin'} (a more general form of Theorem~\ref{Tsumwin})
and Theorem~\ref{Twin}---that will pretty easily give Theorem~\ref{T1/e},
as well as playing a role in Theorem~\ref{TKS'}(c).

Section~\ref{GW} is an interlude proving Theorem~\ref{Tgapw};
this doesn't depend on 
Sections~\ref{Geometry}-\ref{Consequences} but does make use of
well-known correlation inequalities of Shepp and Fishburn.

Sections~\ref{P1/e}-\ref{CandA} prove, respectively, Theorems~\ref{T1/e},
\ref{TKS'},
and \ref{TKS}(d).  Theorem~\ref{TKS'}
is largely based on Theorem~\ref{TkKom}, the aforementioned extension
of the main result of \cite{Komlos}, which is proved in 
Section~\ref{Komlos}.

Finally, Section~\ref{Examples} gives examples justifying two assertions from earlier
sections, and Section~\ref{Conjectures} collects a few of the most interesting unanswered 
questions arising from the present work.

\mn
\textbf{Usage} 

\nin
%\emph{Posets}

As above we use ``$<$'' for order in $P$ and ``$\prec$'' 
(or, if necessary, $\prec_\vs$)
for order in the uniform extension $\gss$,
in each case adding a dot for \emph{cover} relations (e.g.\ $x\recovby y$
means $x\leq w\leq y$ $\Ra $ $w\in \{x,y\}$).
As usual $X<Y$ means $x<y$ $\forall (x,y)\in X\times Y$.
Elements $x,y$ of $P$ are \emph{comparable}, denoted $x\sim y$ if 
$x<y$ or $x>y$, and the \emph{comparability graph}, $G(P)$, is the graph on $P$
with adjacency ``$\sim$.''
The \emph{dual} of $P$, denoted $P^*$, has the same ground set as $P$ but reverses 
the order (i.e.\ $x<_{P^*}y$ $\Lra$ $y<_P x$).

A \emph{subposet} is $Y\sub P$ with the natural order.
Note that $P=X\sqcup Y$ means the \emph{set} $P$
is the disjoint union of $X$ and $Y$,
and allows relations between $X$ and $Y$.
A \emph{chain} is a linearly ordered set,
an \emph{antichain} is a set with no relations $x<y$,
and chain, antichain \emph{of P} have the natural meanings.
We use $\A(P)$ for the set of antichains.
The \emph{height}, $\height(P)$, is the maximum number of elements in a chain of $P$.

The sets of 
maximal and minimal elements of $X\sub P$ are denoted $\max (X)$ and $\min (X)$.
An \emph{ideal} (a.k.a.\ \emph{downset}) of $P$ is $I\sub P $ satisfying
$x<y\in I\Ra x\in I$, and \emph{filter} (\emph{upset}) is defined similarly.
We use $\langle X\rangle $ and $\llangle X \rrangle$ for the ideal and filter
generated by $X$ (e.g.\ $\langle X\rangle =\{x:\exists y\in X, x\leq y\}$), 
and as earlier set $\ga(x)=|\lxr|$ and $\gb(x) =|\llxr|$.

As above, $E(P)$ is the set of linear extensions (we will usually say simply \emph{extensions})
of $P$ and $e(P)$ is its size, and 
we use either $f$ or $\gss$ for a uniform member of $E(P)$ (depending on viewpoint),
writing $\prr(\cdot)$ for probabilities associated with $\gss$.
A central statistic is $h(x):=\E f(x)$, the \emph{average height} of $x$.

\nin
%\emph{General}

As usual $[n]=\{1\dots n\}$, $2^V$ is the power set of $V$ 
(which is ordered by containment), and
we use $\{\cdots\}$ for events.
% and $\sqcup$ for disjoint union.
We will sometimes write $A\sm x$ and $A\cup x$ for $A\sm \{x\}$, and $A\sm \{x\}$.
We will tend to omit $[ ~]$'s in expectations,
e.g.\ using $\E X$ and $\E XY$ rather than $\E[X]$ and $\E[XY]$.
We use standard asymptotic notation (e.g.\ $o(\cdot)$, $\gO(\cdot)$) as 
well as $a \gg b$ for $a=\go(b)$ and $a\more b$ for $a > (1-o(1))b$,
with the parameter underlying the limits (usually $n$ or $w(P)$) given by context.
We usually  
pretend large numbers are integers, rather than clutter the discussion with
irrelevant floor and ceiling symbols.

\section{Parameters}\label{Parameters}

\nin

We would like to put Conjecture~\ref{CKS}
in a broader context, still following the
idea that \eqref{1/2} can fail only for quite special $P$,
and, more generally, hoping
to better understand the mysterious space $E(P)$.
To this end we introduce several new parameters 
associated with a poset.
At least one of these---``win''---is crucial for our balance results,
but we also think the parameters themselves are interesting, and
will devote some effort to
understanding relations between them.
(We postpone discussion of continuous versions, though we 
will later often work with these.)

We use 
$\gs^2(x) $ for the variance of $f(x)$ (and $\gs(x)$ for its standard deviation), and add:
\beq{win}
\win(x) = \E [r(x)-q(x)]
\enq
($\win$ for ``window''), 
with 
\beq{qr}
\mbox{$r(x) = \min\{f(y): y>x)\}$ and $q(x) =\max\{f(y): y<x)\}$},
\enq
where we define the min to be $n+1$ if $x\in \max(P)$ 
and the max to be $0$ if $x\in \min(P)$;
\beq{pi}
\mbox{$\pi(x) =|\Pi(x)|$, where $\Pi(x)=  \{y\in P-x: y\not\sim x\}$};
\enq
\beq{tau}
\tau(P) =\max_{X\sub P} |X|^{-1}H(\gss|_X),
\enq
where $H$ is (say binary) entropy; and
\beq{gap}
\gap(P)=\max\{h(x_1),h(x_2)-h(x_1)\dots h(x_n)-h(x_{n-1}),
n+1-h(x_n)\},
\enq
where the elements of $P$ are $x_1\dots x_n$ with 
$h(x_1)\leq \cdots \leq h(x_n)$.

For $\ga\in\{\gs,\pi,\win\}$, we set
$\ga(P)=\max\{\ga(x):x\in P\}$.

Part of the original motivation here was the idea that each of the 
above parameters
%$\gap,\tau,\win,\gs,\pi$ 
could
replace width in Conjecture~\ref{CKS}, in 
the sense that \eqref{1/2} should hold if any one of them tends to infinity.
The least plausible of these possibilities---which easily implies all the others, including
Conjecture~\ref{CKS}---is an
old 
%, unpublished
conjecture of the second author (from before 1990 if memory serves):
\begin{conj}\label{maxC}
If $~\pi(P)\ra\infty$ then \eqref{1/2} holds.
\end{conj}
\nin
This was recently shown by the first author \cite{Aires}
for $P$ of width 2, and we can now show it for any fixed width
(see Corollary~\ref{Cpi}).
The new argument is different from, and considerably 
easier than, the original;
but it is still not easy and is 
%mostly 
postponed to \cite{AK2}.

%%%%%%%%%%%%%%%%%%%
\iffalse
\mn
As we will see shortly, two of these alternate conjectures (those with $\win$ and $\gs$) 
are equivalent to  Conjecture~\ref{CKS}, and two ($\gap$ and $\tau$)
are weaker.  
\fi %%%%%%%%%%%%%%%%%

For discussing relations between our parameters, with which we now include $w$
(so, to repeat, the parameters are 
$w, \gs, \win,\pi, \tau,$ and $\gap$),
we use $\gc\leadsto\vt$ to mean that $\vt$ tends to infinity if $\gc$ does.
Some of these implications are trivial or nearly so, namely
\beq{easy.arrows}
\mbox{$\tau\leadsto w$, $\,\, \win\leadsto \gs$,}
\enq
and $\gc\leadsto \pi$ for $\gc\in \{\gs,\win,w,\tau,\gap\}$
(``nearly so'' is for $\win\leadsto \gs$; see \eqref{varf});
but the next three
are more substantial.
\begin{thm}\label{Twwin}
$\,\,w\leadsto \win$.
\end{thm}
\begin{thm}\label{Tgapw}
$\,\,\gap \leadsto w$.
\end{thm}
\nin
So we also have $w\leadsto\gs$ and $\gap\leadsto \win$.
\begin{thm}\label{Ttaugs}
$\,\,\pi \leadsto \gs$.
\end{thm}

\nin
On the other hand it's easy to see 
that there are no $\leadsto$'s for the sequence
$(\gs,\win,w,\tau,\gap)$, and that $\pi\not\leadsto \win$ 
(so also $\pi\not\leadsto w,\tau,\gap$).
In sum, 
\beq{heirarchy}
\left.
\begin{array}{rr}
\tau\not \leadsto\gap\\
\tau
\end{array}
\right\}\leadsto
w\leadsto \win\leadsto \gs\leadsto \pi \leadsto \gs,
\enq
leaving just one unresolved possibility:
\begin{conj}\label{Qparams}
$\,\,\gap\leadsto \tau$.
\end{conj}
\nin
\emph{Some} supporting evidence for this is provided by Theorem~\ref{TKS'}(b),
which implies that $\tau\ra\infty$ if, in \eqref{gap}, $h(x_1)\ra \infty$.

Theorems~\ref{Twwin} and \ref{Tgapw} are discussed further at the end
of this section and proved in 
Sections~\ref{Consequences}  and \ref{GW}, with the 
more difficult proof of
Theorem~\ref{Ttaugs} given in \cite{AK2}.
Combined with Theorem~\ref{Twsig} below, Theorem~\ref{Ttaugs} has the following 
surprising consequence (as discussed after the statement of Corollary~\ref{Cpi}).
%Theorem~\ref{Twsig}).
\begin{thm}\label{Tequiv}
Conjecture~\ref{CKS} implies Conjecture~\ref{maxC}.
\end{thm}
\nin
Thus, in view of \eqref{heirarchy}, Conjecture~\ref{CKS} and the seemingly far stronger 
Conjecture~\ref{maxC}, together with the
the $\win$ and $\gs$ versions of Conjecture~\ref{CKS}, are all
equivalent (and imply the
gap and $\tau$ versions).

\begin{thm}\label{Twsig}
For fixed $w$, $P$ of width at most $w$, and $ x\in P$,
\[
\mbox{$\gd_x\ra 1/2$ as $\gs(x)\ra\infty$.}
\]  
In particular, if $\gs(P)\ra\infty$ then \eqref{1/2} holds.
\end{thm}
\nin
The proof of Theorem~\ref{Twsig} will be given, with that of Theorem~\ref{Ttaugs},
in \cite{AK2}.  (The proofs are unrelated, but 
the theorems form a sort of team in combining 
%having the theorems in one place feels natural, since they combine 
to produce Theorem~\ref{Tequiv}.)

\nin

The general idea behind Conjecture~\ref{CKS} was that failure of \eqref{1/2}
for large $P$ can only occur in quite limited situations (and the conjecture says
that large width excludes these).  
Theorem~\ref{Twsig} provides a nice instantiation of this idea in the setting of \emph{bounded} 
width, as does the following, perhaps more striking, consequence in the direction of Conjecture~\ref{maxC}.

\begin{cor}\label{Cpi}
For fixed $w$ and $P$ of width at most $w$, 
if $\pi(P)\ra\infty$ then \eqref{1/2} holds.
\end{cor}

For Theorem~\ref{Tequiv}, just notice that
Conjecture~\ref{CKS} implies its $\gs$
version---which by Theorem~\ref{Ttaugs} implies Conjecture~\ref{maxC}---since Theorem~\ref{Twsig} reduces proving the $\gs$ version to 
proving it for large width.
(Precisely:  for any $\eps>0$, Conjecture~\ref{CKS} gives $\gd(P)> 1/2-\eps$ if
$w(P)>w_\eps$, and Theorem~\ref{Twsig} gives the same conclusion for $P$ with $w(P)\le w_\eps$
and $\gs(P)$ sufficiently large.)
%; and the other equivalences then follow from \eqref{heirarchy}.

\nin
\emph{Further remarks.} 
1.  The implications $w\leadsto \win\leadsto \gs\leadsto \pi$ are ``local,'' the last two
in the natural sense (e.g.\  large $\win(x)$ implies  large $\gs(x)$),
%and similarly  for $\gs$ and $\pi$), 
and the first in that a large antichain must 
contain an element with large $\win$.  On the other hand the
conjectured consequences for \eqref{1/2} are \emph{not} local;
e.g.\ 
a large antichain may contain no well-balanced pair, 
and $\pi(x)$ large does not imply that $x$ belongs to a well-balanced pair.

\nin
2.  Of course proving $\gc\leadsto\vt$ means bounding $\vt$ below by some 
function of $\gc$; but but we will usually leave such functions unspecified except where 
they are nice or needed in proofs.

\mn

We close this section with two crucial points involving win 
(both proved in Section~\ref{Consequences}), and a favorite
conjecture related to Theorem~\ref{Tgapw}.
The first ``crucial point'' is
a quantified version of Theorem~\ref{Twwin},
%which is crucial for Theorem~\ref{TKS'}(c) and most of Theorem~\ref{T1/e}.
which is tight when $P$ is itself an antichain.
\begin{thm}\label{Tsumwin}
For any $P$ and $A\in \A(P)$,
\[   %beq{sumwin}
\mbox{$\sum_{x\in A}\win(x) \geq (n+1)|A|^2/n.$}
\]   %enq
\end{thm}
\nin
In particular, $\win(P)\geq (n+1)w(P)/n$, yielding Theorem~\ref{Twwin}.

The second point, the following connection between win and \eqref{1/e}, is the basic mechanism
underlying Theorem~\ref{T1/e}.
\begin{thm}\label{Twin}
For $X\sub P$, if 
\beq{winsum}
\sum_{x\in X}\win(x)\gg n,
\enq
%---equivalently, if $\sum_{x\in X}\Win(x)\gg 1$---
then 
there are $x,y\in X$ with $\gd_{xy}>1/e-o(1)$.
\end{thm}

\mn

In connection with Theorem~\ref{Tgapw} we mention the 
following much stronger possibility, another old conjecture of the second author.
(This one
%---precisely, just for $|\max(A)|=2$---
\emph{has}
appeared previously, in \cite{Brightwell-Trotter,Biro-Trotter} which 
it partly motivated.)

\begin{conj}\label{A1}
For any ideal A in a poset P, 
\beq{ideal}
\max\{h(x):x\in A\}\geq  |A|-|\max (A)|+1.
\enq
\end{conj} 
\nin
(Note this would imply $\gap(P)\leq 2w(P)-1$, a strong form of Theorem~\ref{Tgapw}.)
Conjecture~\ref{A1}
is true, and best possible, when $A=P$, 
which ought to be the worst case
(its truth is not obvious but follows easily
from \cite{Stanley1} via Lemma~\ref{LSKS}(a) below), 
but in general we don't know that $|A|-\max\{h(x):x\in A\}$
is bounded by \emph{any} function of $|\max(A)|$, 
even when $|\max(A)|=w(P)=2$. 

\nin

While it may seem surprising that the innocent-looking ``$\gap\ra\infty$'' should 
be at all strong (e.g.\ strong enough to imply large width as in Theorem~\ref{Tgapw}), 
we believe it is stronger still:
\begin{conj}\label{CgapkKom}
For fixed $k$, if $\gap(P)\ra \infty$ then \eqref{1/k!} holds.
\end{conj}
\nin
(Note Theorem~\ref{TKS'}(b) is the special case where $h(x_1)$ in \eqref{gap} is large.)
The conjecture clearly fails if we replace gap by 
$\pi, \gs,\win$, or $w$, but could be true with $\tau$, a stronger statement if, 
as in Conjecture~\ref{Qparams},
$\gap\leadsto\tau$.
(Note that Conjecture~\ref{CgapkKom} is in fact a strong form of Conjecture~\ref{Qparams}.
It may also be worth noting that for a \emph{general} random $\gss\in \mS_n$ 
even $H(\gss) > .9n\log n$ doesn't guarantee balance as in \eqref{1/2}; cf.\ 
\cite{Leighton-Moitra,ADK}.)
We also believe that any of the conditions in Theorem~\ref{T1/e} should be enough
for \eqref{1/k!}, though we don't even know they're enough for \eqref{1/2}.

\section{Geometry}\label{Geometry}

We use \emph{body} to mean a full-dimensional compact convex subset of $\mbR^d$
(for some $d$), 
and denote Euclidean volume by $|\cdot|$, with dimension 
given by context.

We will be interested in two polytopes associated with $P$.  The first of these, 
the \emph{order polytope}, is
\[
\Oh(P) = \{v\in [0,1]^P: x< y\Ra v_x\leq v_y\},
\]
and the second, the \emph{chain polytope}, will be discussed shortly.
(See \cite{Stanley2} for background on both.)
We will usually capitalize 
continuous counterparts of discrete quantities, so to begin use 
$F$ for a uniform point of $\Oh(P)$.

%The connection between $\Oh(P)$ and extensions is clear, since the
The extensions of $P$ naturally decompose
$\Oh(P)$ into simplices of volume $1/n!$; namely, $f\in E(P)$, say with $f(x_i)=i$
corresponds to the simplex 
\beq{simplex}
\{0\leq v_{x_1}\leq \cdots\leq v_{x_n}\leq 1\},
\enq
yielding in particular $\pr(F(x)<F(y))= \pr(f(x)<f(y))$ and
\beq{eOhP}
|\Oh(P)| = e(P)/n!
\enq
If we think of choosing $f$ and then $F$ uniformly from the corresponding simplex,
then we see that $\E [F(x_i)|f] = i/(n+1)$ and the centroid of $\Oh(P)$ is
\beq{H}
H(x):=\E F(x) =\E f(x)/(n+1) = h(x)/(n+1).
\enq

We also use
\beq{QR}
\mbox{$Q(x) =\max_{y< x} F(y)~$ and $~R(x)=\min_{y> x} F(y)$}
\enq
(with $Q(x)= 0$ if $x$ is minimal and $R(x)=1$ if $x$ is maximal; cf.\ \eqref{qr}), and
\beq{Win}
\Win(x)=\E[R(x)-Q(x)]=\E[r(x)-q(x)]/(n+1) =\win(x)/(n+1).
\enq

For a given setting of $(F(z): z\neq x)$, $F(x)$ is uniform from $[Q(x),R(x)]$, whence
\beq{HQR}
H(x)=\E[Q(x)+R(x)]/2
\enq
---in particular
\beq{WinH}
\mbox{$H(x)=\Win(x)/2 ~$ if $~x\in \min (P)$,}
\enq
since in this case  $Q(x)\equiv 0$---and
\beq{EFQ}
\E[F(x)-Q(x)]=\E[R(x)-F(x)]=\Win(x)/2.
\enq
Moreover, since $\xi$ chosen uniformly from $[0,a]$ satisfies $\E|\xi-b|\geq a/4$ (for any $b$) and 
$\var(\xi) =a^2/12$, we have
\beq{EFxy}
\E|F(x)-F(y)| \ge \E[R(x)-Q(x)]/4=\Win(x)/4
\enq
(equivalently, $\E|f(x)-f(y)|\geq \win(x)/4$) for any $y\neq x$, and
\beq{VarF}
\Var(F(x))\geq \E\{\Var[F(x)|Q(x), R(x)]\}  \geq \E (R(x)-Q(x))^2/12 
\geq \Win^2(x)/12.
\enq
The discrete version of \eqref{VarF} is 
\beq{varf}
\gs^2(x) \geq (\win(x)-1)^2)/12,
\enq
which justifies
our assertion (in \eqref{easy.arrows}) that $\win\leadsto \gs$.
(The justification is not immediate from \eqref{VarF}, due to the last expression in
\beq{Var}
\Var(F(x)) = \frac{\gs^2(x)}{(n+1)(n+2)} + \frac{H(x)(1-H(x))}{n+2}
\enq
[routine derivation omitted], which roughly corresponds to the residual variance in $F(X)$ given $f(x)$.)

\mn
\emph{Remark.}
Though we don't have a concrete use for it at the moment,
the following counterpart of \eqref{VarF}
seems basic,
but has proved remarkably difficult to attack.

\begin{conj}\label{CWinVar}
There is a fixed $\eps>0$ such that for every $P$ and x, $\Win(x)>\eps\Var(F(x))$.
\end{conj}
\nin
This is best possible (up to the value of $\eps$; 
see also Conjecture~\ref{CWinVar'}) when $P$ is a chain with $x$ not too near the ends.
On the other hand, we don't even know that one can't simultaneously have 
$\Var(F(x)) = \gO(1)$ and $\Win(x)=O(1/n)$ (respectively the 
largest and smallest possible values of these parameters).

%\mn

We now turn to the chain polytope, 
\[
\cee(P) = \{v\in [0,1]^P: \mbox{$\sum_{x\in C}v_x\leq 1$ for every chain $C$ of $P$}\}.
\]
Though first considered as an object of particular interest 
%The chain polytope per se first appears in the literature 
only in \cite{Stanley2}, $\cee(P)$ is
also the case $G=\GP$ of the \emph{vertex-packing polytope} of a graph $G$,
\beq{VPP}
\conv(\{\mbone_{I}: I\in \I(G)\}) \sub [0,1]^{V(G)},
\enq
where $\I(G) $ is the collection of independent sets of $G$
(e.g.\ $\I(\GP) =\A(P)$).

In what follows we will sometimes work with the more general 
notion of a \emph{convex corner} (we will say simply \emph{corner});
that is, a compact convex $\K\sub (\mbR^+)^n$ 
with nonempty interior, satisfying
$\uzero \leq \uv\leq\uw\in \K$ $\Ra$ $\uv\in \K$
(where ``$\leq$'' is product order).
Where valid, we will usually give results for corners, since
this costs nothing and may be helpful in showing what we actually use.
(Of course a reader may just think of $\cee(P)$.)
In all cases the ground set is $V$, with $|V|=n$; $x$, $y$ are in $V$
and $\uv$ is in $\mbR^V$;
and, with $\K$ the body in question,
$\uc$ is the centroid of $\K$ and
$\ud=(d_x:x\in V)$ is given by
\beq{d_x}
d_x = |\K|/|\K-x|,
\enq
where $\K-x$ is the intersection of $\K$ with the hyperplane $\{v_x=0\}$.
%the $x^{\text{th}}$ coordinate hyperplane.

\mn

For posets, a basic observation of Stanley \cite{Stanley2} is that there is a 
bijective, piecewise-linear, volume-preserving 
``transfer map'' $F\mapsto F^*$ from $\Oh(P)$ to
$\cee(P)$ given by $F^*(x)= F(x)-Q(x)$.  This implies
\beq{ecee}
|\cee(P)| =|\Oh(P)| \,\,(= e(P)/n!),
\enq
says that the centroid of $\cee(P)$ is given by 
\beq{centWin}
c_x = \Win(x)/2
\enq
(see \eqref{EFQ}), and gives an alternate expression for $\ud$:
\beq{d.alt}
d_x = |\cee(P)|/|\cee(P-x)| =e(P)/(ne(P-x)).
\enq
As we will see in Lemma~\ref{TApWin},
$\ud$ turns out to be a good proxy for $\uc$, so also for Win if $\K=\cee(P)$.

\section{Log-concavity}  \label{Log-concavity}

The geometric viewpoint gets much of its power from the tractability of 
log-concave functions.  Here we briefly review what we need from 
this subject.  (See also the lucid account in \cite{LV}.)

Recall that $g:\mbR^n\ra\mbR^+$ is \emph{log-concave} if 
\[
g(\gl x +(1-\gl)y)\geq g(x)^\gl g(y)^{1-\gl}
\]
for all $x,y\in \mbR^n$ and $\gl\in[0,1]$.
We call a random $F:\mbR^n\ra\mbR$ \emph{lc-distributed}
if its density function is log-concave.
We also say a random $\mbZ$-valued $f$ is lc-distributed
if the sequence $\{\pr(f=k)\}$ is log-concave; that is, 
$\pr(f=k)^2\geq \pr(f=k-1)\pr(f=k+1)$ for all $k$.

The following version of the Pr\'ekopa-Leindler Inequality,
is contained in Theorem~5.1 of \cite{LV}.

\begin{thm}\label{TDLP}
If $g:\mbR^n\ra\mbR^+$ is log-concave, then its projection on 
any subspace of $\mbR^n$ is log-concave.
\end{thm}
\nin
(We will use only one-dimensional projections, but see \cite{Kahn-Yu}
for a higher-dimensional poset application.)

\mn

For the rest of this section $Z$ is an $\mbR$-valued
r.v.\ with density $g$.  
%We begin with the trivial:
%\begin{prop}\label{Ptriv}
%If $Z\geq 0$ and $g$ is nonincreasing on $\mbR^+$, then
%$\E Z\geq 1/(2g(0))$.
%\end{prop}
The first part of the next lemma is trivial and the second is
contained in \cite[Lemma~5.3(c,d)]{LV}.
\begin{lemma}\label{LLV}
For nonnegative $Z$:

\nin
{\rm (a)} If $g$ is nonincreasing on $\mbR^+$, then
$\E Z\geq 1/(2g(0))$.

\nin
{\rm (b)} 
If g is log-concave then 
$\E Z\leq 1/g(0)$ and $ \E^2 Z \leq 2\E^2Z.$
%$~1/g(0) \geq \E Z\geq \E^2 Z/2.$
\end{lemma}

We briefly depart from geometry 
to 
recall two observations in the same vein but for the discrete setting.
(These are standard, but we don't know a reference.)
\begin{prop}\label{Ldisclc}
If $f$ takes values in the positive integers 
and is lc-distributed, then

\nin
{\rm (a)} $\E f\leq 1/\pr(f=1)$;

\nin
{\rm (b)}  $\Var(f) =O((\max_k\pr(f(x)=k))^{-2})$.
\end{prop}
We return to $Z$.
The next result (perhaps in more precise form) is presumably not new.
\begin{lemma}\label{T1.3}
If Z is lc-distributed then for any $\eps$,
$\pr(|Z|\le \eps\E|Z|/3) \leq \eps$.
\end{lemma}

\nin
\emph{Proof.}
It is enough to show this when $Z$ has compact support and $\E|Z|=1$;
so we assume this,
let $Z$ maximize $\pr(|Z|\le \eps\E|Z|/3)$,
%the probability in question 
 and
use $p$ for the density of $Z$ (so $\log p$ is concave). 
We claim that $\log p$ is linear on each of $(-\infty,-\eps/3]$, $[\eps/3,\infty)$.
For example, if the second linearity fails, then we may replace $p$ by a 
density $q$ that agrees with $p$ on $(-\infty, \eps/3)$, while on $[\eps/3,\infty)$, $\log q$
is linear and $q$ has the same integral as $p$.
Then $Z'$ with density $q$ has $\E|Z'|> \E|Z|=1$ and 
$\pr(|Z'|\le \eps\E|Z'|/3)> \pr(|Z'|\le \eps/3)=\pr(|Z|\le \eps/3)$,
contradicting our choice of $Z$.

Let $p(t_0)=\max_tp(t)$; by the above linearity, $t_0\in[-\eps/3,\eps/3]$. 
If $p(t_0)\le 3/2$, then 

\[
\pr\left(|Z|\le \eps/3\right)\le (3/2)(2\eps/3)= \eps;
\]
and otherwise, Lemma~\ref{LLV}(b) gives another contradiction:
\[
\E|Z|\le |t_0|+\E|Z-t_0|\le |t_0|+1/p(t_0)\le \eps/3+2/3<1.
\]
\qed

%%%%%%%%%%%%%%%%%%
\iffalse

\nin
Then $Z'$ with density $q$ has $\E|Z'|> \E|Z|=1$ and 
$\pr(|Z'|\le \frac{\eps}{3}\E|Z'|)> \pr(|Z'|\le \frac{\eps}{3})=\pr(|Z|\le \frac{\eps}{3})$,
contradicting the choice of $Z$.

Let $p(t_0)=\max_tp(t)$; by the above linearity, $t_0\in[-\frac{\eps}{3},\frac{\eps}{3}]$. 
If $p(t_0)\le \frac{3}{2}$, then 

\[
\pr\left(|Z|\le \frac{\eps}{3}\right)\le \frac{3}{2}\frac{2\eps}{3}= \eps;
\]
and otherwise, Lemma~\ref{LLV}(b) gives another contradiction:
\[
\E|Z|\le |t_0|+\E|Z-t_0|\le |t_0|+\frac{1}{p(t_0)}\le \frac{\eps}{3}+\frac{2}{3}<1.
\]
\fi
%%%%%%%%%%%%%%%%%%

The following fundamental fact is Lemma~5.4 of \cite{LV}.
\begin{thm}\label{TGrunLV}
For any lc-distributed $Z$, $\pr(Z \geq \E Z) \geq 1/e$.
\end{thm}
\nin
In view of Theorem~\ref{TDLP}, this generalizes
the result of Gr\"unbaum mentioned in Section~\ref{Intro},
which we will also use below 
%(in Section~\ref{P1/e}) 
and now state formally:
\begin{thm}\label{TGrunbaum}
For a body $K\sub \mbR^d$ and $H^+$
one side of a hyperplane H through the centroid of K,
\[
|H^+\cap K|\geq |K|/e.
%$H=\{x\in \mbR^d: \langle \ua,x\rangle
\]
\end{thm}

The following variant of Theorem~\ref{TGrunLV} will be crucial (see Lemma~\ref{LWin}).

\nin
\begin{lemma}\label{LZZ}
If $\E |Z|\gg |\E Z| $ then $\pr(Z>0) \more 1/e$.
$ ~~ %\beq{lc1}
%\E |Z|\gg |\E Z| \,\, \Ra \,\,\pr(Z>0) \more 1/e.
$   %\enq
\end{lemma}

\nin
\emph{Proof.} 
We first observe that the hypothesis of the lemma implies $\gs_Z\gg |\E Z|$, since
\[
\gs_Z^2 =\E |Z|^2-\E^2Z\geq \E^2|Z|-\E^2Z \sim \E^2|Z|\gg \E^2Z.
\]
So it's enough to show
\beq{lc3}
\gs_Z \gg |\E Z| \,\, \Ra \,\,\pr(Z>0)\more 1/e.
\enq
Here 
Theorem~\ref{TGrunLV} says we may assume $\mu:=\E Z<0$, and that, with
$I=[\mu,0]$, it's enough to show
\[   %\beq{ZinI}
\pr(Z\in I)\ll 1.
\]   %\enq
To see this, set $\pr(Z\in I)=\gc$ and
let $x_0$ maximize $g$ on $I$ (so $g(x_0)\geq -\gc/\mu$).
With $h$ the restriction of $g$ to $[x_0,\infty)$ 
scaled to integrate to 1,
Lemma~\ref{LLV} gives
\[
E[Z^2|Z\geq x_0] \leq 2\E^2[Z|Z\geq x_0] \leq 2h(x_0)^{-2} \leq 2g(x_0)^{-2}
\leq 2(\mu/\gc)^2,
\]
and the same argument (with a different $h$) gives $E[Z^2|Z\leq x_0]\leq 2(\mu/\gc)^2$.  
So $\gs_Z\leq \sqrt{2}|\mu|/\gc$, 
which with 
$\gs_Z\gg|\mu|$ 
(as in \eqref{lc3})
gives
the desired $\gc =o(1)$.

\qed

\section{Consequences}\label{Consequences}

This Section serves as a bridge between the material of Sections~\ref{Geometry}
and \ref{Log-concavity}, and the balance results of Sections~\ref{P1/e}
and \ref{SB}.
The main points are Theorem~\ref{Tsumwin'} (which contains
Theorem~\ref{Tsumwin} and uses Theorem~\ref{TEHS'}
and Lemma~\ref{TApWin});
Theorem~\ref{Twin} (an easy consequence of Lemma~\ref{LWin},
but the basis for Theorem~\ref{T1/e}); and 
Lemma~\ref{Lcentroid}, a corollary of Theorem~\ref{TGrunbaum} that 
will be needed for Theorem~\ref{T1/e}(d).

To begin we note three poset consequences of results of Section~\ref{Log-concavity}
(namely, Theorem~\ref{TDLP} contains (a), and, given (a), Theorem~\ref{TGrunLV} implies (b) 
and Lemma~\ref{T1.3}, with \eqref{EFxy}, implies (c)):
\begin{cor}\label{for.posets}
For any $x,y\in P$,

\nin
{\rm (a)}
$F(x)-F(y)$ is lc-distributed;

\nin
{\rm (b)}
if $h(x)\leq h(y)$ then $\pr (f(x)<f(y))\geq 1/e$;

\nin
{\rm (c)}
$\pr(|F(x)-F(y)|\leq \eps\Win(x)/12) < \eps$.
\end{cor}

The next two inequalities follow from
Proposition~\ref{Ldisclc} and
log-concavity results of \cite{Stanley1,Kahn-Saks},
and will be used for Theorem~\ref{TKS'}(b) and Theorem~\ref{Tgapw}
respectively.  These are surely not new, but we don't know a reference so 
give the easy proofs.

\begin{lemma}\label{LSKS}
{\rm (a)}  For $x\in P$, $h(x)\leq 1/\pr(f(x)=1)$.

\nin
{\rm (b)}  For distinct $x,y\in P$, $h(y)-h(x)\leq 1/\pr(f(y)-f(x)=1)$.
\end{lemma}

\nin
\emph{Proof.}
{\rm (a)}  This is contained in Proposition~\ref{Ldisclc}, since
log concavity of $\{\pr(f(x)=k)\}$ was shown in \cite{Stanley1}.

\nin
{\rm (b)} Here the relevant result, from \cite{Kahn-Saks}, says that
$\{\pr(f(y)-f(x)\}_{k\geq 1}$ is log-concave.  
Combined with Proposition~\ref{Ldisclc} this gives, with $B=\{f(y)>f(x)\}$,
\begin{eqnarray*}
h(y)-h(x)&\leq  &\pr(B)\E [f(y)-f(x)|B] \leq  \pr(B)/\pr(f(y)-f(x)=1|B)\\
&=&
\pr(B)^2/\pr(f(y)-f(x)=1) \leq  1/\pr(f(y)-f(x)=1).
\end{eqnarray*}
\qed

The following observation for a corner $\K$ gives the relation between $\uc$ (the centroid)
and $\ud$ 
(see \eqref{d_x})
alluded to at the end of Section~\ref{Geometry}.
Here we are mainly interested in the lower bound, which, as said above,
will be used for Theorem~\ref{Tsumwin}.

\begin{lemma}\label{TApWin}
For every convex corner $\K$ and $x\in V$,
\[
d_x/2\le c_x\le d_x.
\]
In particular, when $\K=\cee(P)$, this says
$d_x\leq \Win(x) \leq 2d_x$.
\end{lemma}
\nin
\emph{Proof.}
Projecting $\K$ (i.e.\ uniform distribution on $\K$)
on the $x$-axis produces a nonincreasing log-concave density $g$ with 
mean $c_x$ and $g(0) = |\K-x|/|\K| =1/d_x$.
The first part of the lemma now follows from Lemma~\ref{LLV},
and for the second part we recall from \eqref{centWin} that
$\cee(P)$ has centroid $\Win/2$.

\qed

\nin
\emph{Proof of Theorem~\ref{Tsumwin}.}
We prove this at the level of corners:
\begin{thm}\label{Tsumwin'}
For a corner $\K$ and $I\sub V$ with $\mbone_I\in \K$,
$   %beq{sumwin}
\mbox{$\,\sum_{x\in I}c_x \geq |I|^2/(2n)$.}
$   %enq
\end{thm}
\nin
(This contains Theorem~\ref{Tsumwin} since $\cee(P)$ has centroid $\Win/2$.)

\nin

Theorem~\ref{Tsumwin'} is an easy consequence of the next
inequality, which, as observed in \cite{BB}, 
is implicit in~\cite{Meyer} (see \cite[p.\ 425]{BB}).
\begin{thm}\label{TEHS'}
For any convex corner $\K$ and $\uv\in \K$,
\beq{Meyer}
\mbox{$|\K|\geq n^{-1}\sum_{x\in V} v_x|\K-x|.$}
\enq
\end{thm}
\nin
When $\K$, $I$ are as in Theorem~\ref{Tsumwin'} and 
$v=\mbone_I$, \eqref{Meyer}
becomes
\beq{Meyer'}
\mbox{$|\K|\geq n^{-1}\sum_{x\in I} |\K-x|,$}
\enq
and further specializing to $\cee(P)$
yields the following result of \cite{EHS,Stachowiak},
which with the present argument is enough for Theorem~\ref{Tsumwin}.
\begin{thm}\label{TEHS}
For any $P$ and $A\in \A(P)$,
$
e(P)\geq \sum_{x\in A}e(P-x).
$
\end{thm}
\mn

Combining \eqref{Meyer'} and
the lower bound in Lemma~\ref{TApWin} 
now gives Theorem~\ref{Tsumwin'}: with sums over $x\in I$,
\[
2\sum c_x \geq \sum d_x = \sum |\K|/|\K-x|\geq
|I|^2/n
\]
(since
$\sum (1/z_i)\geq a^2/\sum z_i$ for nonnegative $z_1\dots z_a$).

\qed

Lemma~\ref{LZZ} (with $Z=F(x)-F(y)$; see Corollary~\ref{for.posets}(a))
gives the next observation, which, in easily implying Theorem~\ref{Twin}, is
one of the main points underlying the present work.
\begin{lemma}\label{LWin}
If $~\E|F(x)-F(y)|\gg |H(x)-H(y)|$ then $\gd_{xy} \more 1/e$.
In particular $\Win(x)\gg |H(x)-H(y)|$ implies $\gd_{xy} \more 1/e$
(equivalently, $\win(x)\gg |h(x)-h(y)| \Ra \gd_{xy} \more 1/e$).
\end{lemma}
\nin
(See \eqref{EFxy} for ``in particular.'')

\nin
\emph{Remark.}
Here one may compare \cite{Kahn-Linial} and \cite{Friedman}, which adapt the approach of 
\cite{Grunbaum} to show (respectively)
$|h(x)-h(y)|< 1 \Ra \gd_{xy} > 1/(2e)$ and
$|h(x)-h(y)|\ll 1 \Ra\gd_{xy}\more 1/e$.
Getting balance from small height differences is also the basis for \cite{Kahn-Saks} and
\cite{BFT}.
See Example~\ref{MaxT2.3} for the perhaps surprising fact that even 
a \emph{large} set of elements of \emph{equal} height needn't contain a pair with
balance better than about $1/e$.

\mn
\emph{Proof of Theorem~\ref{Twin}.}
We write Lemma~\ref{LWin}  in the form
\beq{Cform}
\mbox{for each $\eps>0$ there is a $C_\eps\geq 1$ so that
$\Win(x)\geq C_\eps |H(x)-H(y)| $ implies $ \gd_{xy} > 1/e-\eps$,}
\enq
and conclude that if
$\sum_{x\in X}\Win(x)\geq C_\eps$ then there are $x,y\in X$
with $\gd_{xy}>1/e-\eps$.
To see this just observe that the sum of the lengths of the
intervals
\[
I_x:=[H(x)-\Win(x)/(2C_\eps), H(x)+\Win(x)/(2C_\eps)] 
\]
is at least 1, and that $I_x\cap I_y\neq\0$ implies 
$(\Win(x)+\Win(y))/2\geq C_\eps|H(x)-H(y)|$,
which by \eqref{Cform} gives $\gd_{xy}>1/e-\eps$.
(The restriction to $C_\eps\geq 1$ in \eqref{Cform} keeps the $I_x$'s
inside $[0,1]$, since by
\eqref{Win} and \eqref{HQR}, $\Win(x)/2\leq H(x)\leq 1-\Win(x)/2$.)

\qed

\nin
\emph{Remark.}
A similar argument gives Theorem~\ref{TKS}(d) with the target
\eqref{1/2} relaxed to \eqref{1/e}, as follows.
If $C=\{y_1<\cdots < y_k\}$ and
$\gd(P)< 1/e-\eps$, then
the intervals $I_x :=(H(x)-\Win(x)/C_\eps,H(x)+\Win(x)/C_\eps)$ ($x\in A$)
contain no $H(y_i)$'s; so they and the intervals
$(H(y_{i-1},H(y_i))$ comprise 
a disjoint collection of total length at least $|C|-1+\sum |I_x|= n-|A|+1+\sum |I_x| > n$
(where, since at least $|A|-2$ of the $I_x$'s are contained in $[0,1]$, 
the sum is much larger than $|A|$ by Theorem~\ref{Tsumwin}).

\mn

Finally, Theorem~\ref{T1/e}(d) will follow easily from Theorem~\ref{Twin} 
and the next observation.

\begin{lemma}\label{Lcentroid}
For any convex corner $\K$, 
\beq{|K|}
\mbox{$\prod d_x \leq |\K| \leq e\prod (nc_x)/n!  \,\,\,\,\,(< e^n\prod c_x).$}
\enq
\end{lemma}
\nin
(Here 
the upper bound is the main point.
Similar bounds based on graph entropy are given in
\cite{EandS};
see Remark 2 near the end of Section~\ref{P1/e}.)

\mn
\emph{Proof.}
For the upper bound, with
$H$ the hyperplane through $\uc$ orthogonal to 
$\ub := (1/(nc_1)\dots 1/(nc_n))$, Theorem~\ref{TGrunbaum},
gives 
\[    %begin{eqnarray*}
\mbox{$|\K|/e \leq |\{\uv\geq \uzero: \langle \uv,\ub\rangle \leq 1 \}|
= \prod (nc_x)/n!$}
%< e^n\prod c_x$
\]

As to the lower bound, notice that 
$|\K|\geq \prod c_x$ is trivial
(since $\K$ contains $\uc$ and is a corner).  The slight improvement in \eqref{|K|}
follows from the Loomis-Whitney Inequality \cite{Loomis-Whitney},
which for corners says
\[
\mbox{$|\K| \leq \prod |\K-x|^{1/(n-1)}$;}
\]
or, after a little rearranging,
$
\mbox{$|\K|\geq \prod |\K|/|\K-x| = \prod d_x.$}
$
\qed

%\section{Proof of Theorem~\ref{Tgapw}}
\section{Correlation}
%\section{gap $\leadsto$ $w$}
\label{GW}

%The reader may think of this section as a sort of palate cleanser
%preceding the discussion of balance results in Sections~\ref{P1/e} and \ref{SB}.

Before turning to balance we pause for the proof of 
Theorem~\ref{Tgapw}.
% and \ref{Twsig}.
This doesn't need the material developed above but 
does depend on
%Before proving Theorem~\ref{Tgapw} we recall 
two well-known correlation assertions, which we now recall.
The first is Shepp's ``XYZ Theorem'' \cite{SheppII},
which had been a conjecture of Rival and Sands (see \cite{Banff}) and
earlier a question of Kislitsyn \cite{Kislitsyn}.
(The proof of Theorem~\ref{Twsig} in \cite{AK2} again depends on
Shepp's result.)
\begin{thm}\label{TXYZ}
For $x\in P$ and $Y\sub P- x$,
\[
\prr(x\succ Y)\geq \prod_{y\in Y}\prr(x\succ y).
\]
\end{thm}
\nin
(The statement in \cite{SheppII} is for $|Y|=2$, but Theorem~\ref{TXYZ} 
follows via a straightforward induction.
See also \cite{FishburnXYZ} for the surprisingly tricky proof that
equality holds only in trivial cases.)

The second assertion is known as \emph{Fishburn's Inequality} (\cite{FishburnXYZ} or
e.g.\ \cite[Theorem 7.4]{Chan-Pak}), but can also be seen to be contained in
the main result of \cite{SheppI}.
\begin{thm}\label{TFishburn}
If $K$ and $L$ are filters of $P$ then
\[
\frac{e(K\cup L)e(K\cap L)}{e(K)e(L)}\geq \frac{|K\cup L|!|K\cap L|!}{|K|!|L|!}.
\]
\end{thm}

\nin
\emph{Proof of Theorem~\ref{Tgapw}.} 
We prove the following more general statement.
\begin{claim}\label{Cgapw}
If $P=D \sqcup U$ with $D$ an ideal (and $U$ a filter),
$A=\max(D)$,
$B=\min(U)$, and each of $|A|, |B|$ at most $w$, then
\beq{minh(b)}
\min_{b\in B}h(b)-\max_{a\in A}h(a) = O(w^2e^w).
\enq
\end{claim}
\nin
(This of course falls well short of the bound $2w-1$ that
would follow from Conjecture~\ref{A1}.)

\nin
\emph{Proof.}
Here $a$ and $b$ are always assumed to be elements of $A$ and $B$.
Assuming \eqref{minh(b)} fails (or just that $h(b)\geq h(a)$ $\forall a,b$),
%Theorem~\ref{TGrunLV}, with \eqref{FxFylcand} whence,
Corollary~\ref{for.posets}(b) gives
$p(b\succ a) \geq e^{-1}$ $\forall a,b$; so by Theorem~\ref{TXYZ},
\beq{prbsucc}
p(b\succ A) \geq e^{-|A|} \geq e^{-w} \,\,\,\,\, \forall b.
\enq

According to Lemma~\ref{LSKS}(b) it is enough to show that there are $a$ and $b$
with 
\beq{acovbyb}
p(b\cov a) =\gO(w^{-2}e^{-w}).
\enq
Choose $b$ with 
\beq{chooseb}
p(b\prec B\sm \{b\}|A\prec B)\geq 1/|B|
\enq
(e.g.\ choose $b$ to maximize the l.h.s.), and then $a$ with
\beq{prasucc}
p(a\succ A\sm \{a\}|A\prec b) \geq 1/|A|.
\enq
We will show 
\beq{useFish}
p(b\cov a|A\sm \{a\} \prec a\prec b)
\geq 1/|B|,
\enq
which with \eqref{prbsucc} and \eqref{prasucc}
gives \eqref{acovbyb}, since then
\[
p(b\cov a) \geq p(A\prec b)\pr(a\succ A\sm \{a\}|A\prec b)
p(b\cov a| A\sm \{a\} \prec a\prec b) \geq e^{-w}|A|^{-1}|B|^{-1}.
\]

\nin
\emph{Proof of \eqref{useFish}.}
Set $E=\{A\sm \{a\} \prec a\prec b\}$, and
let $K$ run over filters of $U$ containing $b$, setting
%write
%$\pr(\cdot|E,K)$ for $\pr(\cdot|E\wedge \{K=\{x\in U:x\succ a\}\})$.
\[
\psi(K) =p(b\cov a|E\wedge \{\{x\in U:x\succ a\}=K\}).
\]
Then 
\beq{AaBAB}
\psi(U) = p(b\cov a| A\sm \{a\} \prec a\prec B) = p(b\prec B\sm \{b\}|A\prec B) 
\enq
(note $\gss|_A$ and $\gss|_B $ are independent under either of
the conditioning events in \eqref{AaBAB});
so in view of \eqref{chooseb}, 
\eqref{useFish} will follow if we show $\psi(K)\geq \psi(U)$ $ \forall K$.
But 
$
\psi(K) $ (which depends only on the subposet $K$ of $P$) is just
$e(K\sm b)/e(K)$;
so $\psi(K)\geq \psi(U)$ is 
%we want $e(K\sm b)/e(K)\geq e(U\sm b)/e(U)$, i.e.
\[
\frac{e(K\sm b)e(U)}{e(K)e(U\sm b)}\geq 1,
\]
which, since $K \cup (U\sm b)=U$ and $K\cap (U\sm b) =K\sm b$,
is given by Theorem~\ref{TFishburn}.

\qed

\section{Weak balance}
%Proof of Theorem~\ref{T1/e}
\label{P1/e}

Here we prove Theorem~\ref{T1/e}, in which: part
(e) follows 
easily from what we already know (specifically, Theorems~\ref{Tsumwin} and \ref{Twin})
and easily implies
(a) and (b); and for (c) and (d)  
we again use Theorem~\ref{Twin}, its
hypothesis given by 
Claim~\ref{Cl1} for (c),
%(the one point here requiring any effort), 
and by
Theorem~\ref{TGrunbaum} and Lemma~\ref{Lcentroid} for (d).

For (e),  we have (with $A$ and $x$
running over $\A(P)$ and $P$, and the second ``$\geq$'' given by Theorem~\ref{Tsumwin})
\beq{edisp}
\sum\win(x) \geq \sum\win(x)\sum_{A\ni x}\gl_A 
=\sum \gl_A \sum_{x\in A}\win(x) \geq \sum \gl_A|A|^2 \gg n,
\enq
so Thereom~\ref{Twin} gives \eqref{1/e}.

To get (a) from (e), let $\gl$ assign weight 1 to a single antichain of size $\go(\sqrt{n})$.
For (b), partition $P$ into antichains $A_1\dots A_m$, with 
$x\in A_i$ if $i$ is the largest $j$ for which there is a chain $x_1< \cdots< x_j=x$
(so $m=\height(P)$),
and let $\gl$ assign weight 1 to each $A_i$, yielding
\[
\sum \gl_A|A|^2 =\sum |A_i|^2\geq m^{-1}\left(\sum |A_i|\right)^2 = n^2/m\gg n.
\]

\nin
\emph{Remark.}
The above discussion is easily adapted to give the ``local'' versions mentioned in
the remark following
Theorem~\ref{T1/e}:  just restrict the ranges of $A$ and $x$ in \eqref{edisp} 
to $\supp(\gl)$ and $\cup\{A:A\in\supp(\gl)\}$ instead of $\A(P)$ and $P$
(and apply Thereom~\ref{Twin} with $X=\cup\{A:A\in\supp(\gl)\}$).

As noted above, part (c) of Theorem~\ref{T1/e} 
will follow from Theorem~\ref{Twin} once we show
\begin{claim}\label{Cl1}
%There is a fixed $c>0$ such that for very P and x, 
For every P and x, 
$
\,\win(x) \geq 2h(x)/\ga(x).
$
Equivalently, with $\uc$ the centroid of $\cee(P)$,
\beq{cHa}
c(x) \geq H(x)/\ga(x).
\enq
\end{claim}
\nin
(For the equivalence, recall from \eqref{centWin} that $c(x)=\Win (x)/2$.)

\nin
\emph{Proof of \eqref{cHa}}.
It is enough to show \eqref{cHa} given any value of the restriction of $\prec$ to $\lxr$;
so we may assume $\lxr$ is a chain $x_1<\cdots <x_r=x$, with $r=\ga(x)$.
But in this case,
for each $i\in [r-1]$,
\[    %\begin{eqnarray*}
c(x_{i+1})  = H(x_{i+1})-\E Q(x_{i+1}) =H(x_{i+1})-H(x_i)
\geq \E R(x_i)-H(x_i) =c(x_i),
\]   %\end{eqnarray*}
which with $H(x) =\sum c(x_i)$ gives \eqref{cHa}.
\qed

Finally, for (d), using
Lemma~\ref{Lcentroid}
(and again recalling \eqref{centWin}), we have
\[
e(P) ~=~ |\cee(P)|n!
~\leq~ \prod (nc_x) ~=~ (n/2)^n\prod\Win(x) 
%= \left(\frac{n}{2(n+1)}\right)^n\prod \win (x)
~<~ 2^{-n} \prod\win(x);
\]
so $e(P)^{1/n}\gg 1$ (the assumption in (d) implies 
\beq{winwin}
1\ll \prod\win(x)^{1/n} \leq \sum\win(x)/n,
\enq
and Theorem~\ref{Twin} again completes the proof.

\qed

\nin
\emph{Remarks.} 
$1\ll n\prod c_x^{1/n} \leq \sum c_x$.  
1. As mentioned following Lemma~\ref{Lcentroid}, bounds similar to those in \eqref{|K|}
were given in \cite{EandS}.  
Recall that 
the \emph{entropy} of a graph $G$ is
\beq{H(G)}
H(G)= \min\{ n^{-1}\sum_{x\in V}\log (1/a_x): a\in \V(G)\},
\enq
where $\log =\log_2$ and $\V(G)$ is as in \eqref{VPP}.  (Strictly speaking, this is entropy 
w.r.t.\ uniform measure on $V$.
See \cite{Korner} for the origin and \cite{Simonyi} for a good overview.)  
With $G=G(P)$ and $a$ the (unique) minimizer in \eqref{H(G)},
the observation from \cite{EandS} corresponding to \eqref{|K|} is
%says that for $G=G(P)$,
\beq{ent.bds}
\mbox{$n!\prod a_x \leq e(P) \leq n^n\prod a_x $.}
\enq

\nin 
2.  It is also possible to show---this was our original, less easy proof of 
Theorem~\ref{T1/e}(d)---that the condition in (d) implies existence of 
$\gl$ as in (e).  More generally, this is true if 
$\K= \conv(\{\mbone_{I}: I\in \I\}$ with $\I$ an ideal of $2^V$,
where it says:
if $|\K|^{1/n}\gg 1$, then there is a fractional matching $\gl$ of $\I$
satisfying \eqref{gd.frac.m} (with $\A$ replaced by $\I$).

\section{Strong balance}\label{SB}

Here we prove %Theorem~\ref{TKS}
Theorem~\ref{TKS'}(a)-(c).
%strengthened as in .
For this strengthening of the corresponding parts of Theorem~\ref{TKS}
we need the following general principle, of which 
\cite[Thm.\ 1]{Komlos} is the case $k=2$.

\begin{thm}\label{TkKom}
For $k,K$ positive integers and $\eps>0$ there is an $M$ such that: for any
$[k]$-valued r.v.s $X_1\dots X_M$,
there are distinct $i_1\dots i_k\in [M]$
so that for any $\pi\in \mS_{[k]}$ (the symmetric group on $[k]$),
%two permutations $\pi$ and $\gl$ of $[k]$,
\[
|\pr(X_{i_{\pi(1)}}<\cdots <X_{i_{\pi(k)}}) -\pr(X_{i_1}<\cdots <X_{i_k})| <\eps.
\]
So if also $\pr(X_i=X_j)< \eps $ for all distinct $i$ and $j$, then for any 
$\pi\in \mS_{[k]}$,
\[
\pr(X_{i_{\pi(1)}}<\cdots <X_{i_{\pi(k)}}) > 1/k! - C_k\eps.
\]
\end{thm}
\nin
To keep the focus here on Theorem~\ref{TKS'} we postpone the proof of 
Theorem~\ref{TkKom} to Section~\ref{Komlos}.
For the rest of this section we fix $k$, on which our constants will thus depend.

We need the following
immediate consequence of Theorem~\ref{TkKom}.
\begin{cor}\label{CkKom}
Suppose $I_1\cup\cdots \cup I_K$ partitions $[0,1]$ into intervals,
let $X\sub P$, and define r.v.'s $\xi_x$ ($x\in X$) by
setting $\xi_x =\ell$ if $F(x)\in I_\ell$.
If $|X|> M_{K,\eps}$ and 
\[
\mbox{$\pr(\xi_x=\xi_y)< \eps$ for all distinct $x,y\in X$,}
\] 
then $\gd_k(X)> 1/k! - O(\eps)$.
\end{cor}
\nin
(Here and similarly in Corollary~\ref{Lminbalance}, ``If $|X|>M_{K,\eps}$''
is short for ``$\exists M=M_{K,\eps}$ such that if $|X| > M$.''
For ``interval'' and ``partition'' we ignore endpoints.)

A further specialization, 
for use in (a) and (b),
is:

\begin{cor}\label{Lminbalance}
If $X\sub \min(P)$, $|X|> M_\eps$ and, for some $H^*$,
\[
H^*\leq H(x)\leq 2H^*\,\,\forall x\in X,
\]
then $\gd_k(X)> 1/k! -O(\eps)$.
\end{cor}
\nin

\nin
\emph{Proof.}
Let $\nu=\min\{4H^*/\eps,1\}$.  Partition $[0,\nu]$ 
into intervals $I_1\dots I_K$ 
with $|I_j|\leq \eps H^*/12$  ($=\Win(x)/24$; see \eqref{WinH})
and $K=\lceil 12\nu/(\eps H^*)\rceil \leq \lceil 48/\eps^2\rceil$
(a function of $\eps$), and let $I_{K+1} =  [\nu, 1]$.  
Then for $\xi_x$'s as in Corollary~\ref{CkKom}
and distinct $x,y\in X$,
\[
\pr(\xi_x=\xi_y) ~\leq ~\pr(F(x)\geq 4H^*/\eps) + \pr(|F(x)-F(y)|\leq \eps H^*/12)  
~<~ \eps/2+\eps/2~ =~\eps,
\]
with the first $\eps/2$ given by Markov's Inequality
and the second by Corollary~\ref{for.posets}(c).
The lemma is then given by Corollary~\ref{CkKom}.

\qed

\nin
\emph{Proofs of parts (a) and (b) of Theorem~\ref{TKS'}.}
Let $X=\min(P)$.

For (a), noting that $H(x)\geq 1/(n+1)$ $\forall x$
(see \eqref{H}), let $H^*_i =2^i/(n+1)$ for $0\leq i\leq \lceil\log_2 (n+1)\rceil$
and $X_i=\{x\in X: H(x)\in [H^*_i, 2H^*_i)\}$.
Then $\max|X_i|\gg 1$, which implies 
$\gd_k(X) > 1/k!-o(1)$ by Corollary~\ref{Lminbalance}.

For (b) we use Lemma~\ref{LSKS}(a)
(which says $\pr(f(x)=1) \leq 1/h(x)$).
Suppose $\min H(x) = \go/(n+1)$ (with $\go\gg 1$), 
let $H^*_i =2^i\go/(n+1)$ (with $i$ again starting at 0), and define $X_i$ as above.
Then
\[
1= \sum_{x\in X}\pr(f(x)=1) \leq \sum_i|X_i|/(2^i \go)
\]
implies $\max_i|X_i| \gg 1$, and we again finish using Corollary~\ref{Lminbalance}.

\qed

\nin
\emph{Remark.}
A related argument
gives \eqref{1/2} whenever some set of \emph{minimals} 
witnesses large $\tau$ ($\tau$ as in\eqref{tau}):
\begin{thm}\label{Tmintau}
If $\max\{ |X|^{-1}H(\gss|_X):X\sub \min(P)\} \gg 1$, then $P$ satisfies \eqref{1/2}.
\end{thm}

\nin
\emph{Proof of Theorem~\ref{TKS'}(c).}
Fix $\gc_0>0$ and let $A$ be an antichain of size $a = \gc (n+1)$ with $\gc\geq \gc_0$, and
$X =\{x\in A:\win(x) > a/2\}$.
Then $|X| \geq a/2$ (since $|A\sm X|< a/2$ 
by Theorem~\ref{Tsumwin}) and (by \eqref{Win})
\beq{Winwinagain}
\Win(x)>
%\geq a/(2n)=
\gc/2 \,\,\,\,\forall x\in X.
\enq
Given any fixed $\eps>0$, equipartition $[0,1]$ into $\lceil 24/(\eps\gc)\rceil$ intervals 
and define $\xi_x$'s as in Corollary~\ref{Lminbalance}.
Then Corollary~\ref{for.posets}(c) gives, for distinct $x,y\in X$, 
\[
\pr(\xi_x=\xi_y) ~\leq ~\pr(|F(x)-F(y)|\leq \eps\gc/24) < \eps,
\]
ansod Corollary~\ref{CkKom} gives $\gd_k(X)> 1/k! -O(\eps)$
provided ($|X|>$) $\gc n/2> M_\eps$.

\qed

\section{Antichain and chain}
%{Proof of Theorem~\ref{TKS}(d)}
\label{CandA} 

Here we prove Theorem~\ref{TKS}(d), 
aiming for clarity without excessive formality.
To restate:  we are given $P=C\sqcup A$, with $C$ a chain and $A$ an antichain,
and a fixed $\eps>0$,
and should show:
\beq{CAtoshow}
\mbox{there is a $K=K_\eps $ such that $\gd(P)>1/2-\eps$ if $|A|>K$.}
%\mbox{for each $\eps > 0$ there is a $K $ such that $\gd(P)>1/2-\eps$ if $|A|>K$.}
\enq
Throughout this discussion $x, x'$ are in $ A$ and 
$y,z$ (usually with subscripts) are in $C$.
We may assume there are $\hO,\hone\in C$ with
$\hO\le w\le  \hat{1}$ $\forall w\in P$ (this is unnecessary but convenient at \eqref{xinsome}).

Assume for a contradiction that 
\beq{gdPleq}
\gd(P)\leq 1/2-\eps.
\enq
Let 
\[
C=\{\hO= y_k<\cdots <y_0=z_0<\cdots <z_\ell=\hone\},
\]
with $\max (\langle A\rangle\cap C)=\{y_0\}$, and notice that 
$y_0\not > x\not > z_1$ $\forall x$.
We think of ``cells'' $I_s = (y_{s-1},y_s)$ and $J_s=(z_{s-1},z_s)$ ($s\geq 1$),
using e.g.\ ``$x\in I_s$'' for ``$y_s\prec x\prec y_{s-1}$.''
Notice that
\eqref{gdPleq} implies
\beq{xinsome}
\mbox{for each $x$ there is some $I_s$ or $J_t$ that contains $x$
with probability at least $\eps$.}
\enq
(Otherwise, if $w'<w$ are consecutive in $C$ with $\prr(x\succ w')\geq 1/2 > \prr(x\succ w)$,
then $\gd_{x,w} = p(x\succ w)> 1/2-\eps$.)
For each $x$, set
\[
\mbox{$i(x) = \min\{i:x>y_i\}~ $ and $~j(x) = \min\{j:x<z_j\}$.}
\]
\begin{obs}\label{obs1}
For any x, each of $\prr(x\in I_s)$, $\prr(x\in J_s)$ is nonincreasing in $s$.
\end{obs}
\nin
\emph{Proof} (sketch).
We show this for $\prr(x\in J_s)$.
Let $Q=P \cup\{z_{s-1}<x<z_s\}$ (i.e.\ add these relations and close transitively)
and $Q'=P \cup\{z_s<x<z_{s+1}\}$.
Each of these is the disjoint union of $A\sm x$ and a chain with elements $C\cup x$,
say 
$\{\hO= w_1<\cdots< w_m=\hone\}$ for $Q$ and 
$\{\hO= w'_1<\cdots< w'_m=\hone\}$ for $Q'$.

For $x'\in A\sm x$ let
\[
\vp(x') = \max\{i: w_i< x'\}, \,\,  \psi(x') = \min\{i: w_i> x'\};
\]
let $\vp'(x')$, $\psi(x')$ be the corresponding indices for $Q'$;
and \emph{check}:  $\vp'(x')=\vp(x')$ and
\[
\psi'(x')=\left\{\begin{array}{ll}
\psi(x')-1&\mbox{if $j(x')=s$},\\
\psi(x')&\mbox{otherwise.}
\end{array}\right.
\]
So $Q'$ is (isomorphic to) $Q$ with (possibly) some additional relations, and
the observation follows.

\qed

\begin{obs}\label{obs2}
For distinct $x $ and $x'$, either
\[
\max\{i(x),i(x')\} - \min\{i(x),i(x')\} \geq \eps \max\{i(x),i(x')\}
\]
or
\[
\max\{j(x),j(x')\} - \min\{j(x),j(x')\} \geq \eps \max\{j(x),j(x')\}.
\]
\end{obs}
\nin
\emph{Proof}.
If this fails then, with $i =\min\{i(x),i(x')\}$ and
$j =\min\{j(x),j(x')\}$, Observation~\ref{obs1} gives
$
\prr(y_i\prec x,x'\prec z_j) > 1-2\eps.
$
But then 
\[
\prr(x\prec x')~ \geq ~\prr(y_i\prec x,x'\prec z_j)/2 ~> ~1/2 - \eps
\]
and similarly for $\prr(x'\prec x)$, yielding the contradiction
$\gd_{x,x'}> 1/2-\eps$.

\qed

%\mn

Let $a=\max_xi(x)$ and $b=\max_x j(x)$.
Set 
%(pretending large numbers are integers)
$k=\log_{1-\eps} (1/\eps a)$ and partition $[a]$ into
intervals $((1-\eps)^ia,(1-\eps)^{i-1}a]$ for $i\in [k]$ and singletons $\{u\}$ for $u\in [1/\eps]$.
Partition $[b]$ similarly (with $k$ suitably modified) and partition $[a]\times [b]$ into the corresponding
``boxes.''
By Observation~\ref{obs2} no box contains more than one pair $(i(x), j(x))$, so
\beq{xsize}
|A|\leq (\log_{1-\eps} (1/\eps a) +1/\eps)\cdot (\log_{1-\eps} (1/\eps b) +1/\eps).
\enq

Now choose $x$ with $s=s(x):=i(x)+j(x)$ maximum.
By \eqref{xinsome} there is a cell $I$ with $\prr(x\in I)\geq \eps$.
(We may, but needn't, assume, using Observation~\ref{obs1}, that $I\in \{I_1,J_1\}$.)
But 
\[
\prr(x\in I)\leq |A|/s
\]
(since wherever  in $C$ our uniform extension
$\gss$ puts the elements of $A\sm x$, there are at least
$s$ places to insert $x$, at most $|A|$ of which are in $I$).
So $|A| >\eps s$, which, since $s\geq \max\{a,b\}$
(this is all we need from $s$), contradicts \eqref{xsize} for large 
enough $A$.

\qed

\section{A more general Koml\'os Theorem}\label{Komlos}

This will be brief since the main ideas are already in \cite{Komlos}.

\nin
\emph{Proof of Theorem~\ref{TkKom}.}
To avoid pointless precision, we fix $\eta>0$ significantly smaller than $\eps^2$
and use $a\approx b$ for $|a-b| = O(\eta)$ (with the implied constant independent of $\eps$).

An easy application of Ramsey's Theorem allows us to assume,
with $M$ as large as needed, that
%with $\uj$ ranging over $[K]^k$, 
\[
\mbox{\emph{there are numbers $r(\uj)$ for $\uj=(j_1\dots j_\ell)\in [K]^k$ such that
for all $i_1< \cdots < i_k$ in $[M]$,}}
\]
\beq{r(j)}
\pr(X_{i_\ell} = j_\ell \,\, \forall \ell\in [k]) \approx r(\uj).
\enq
It is then enough to show that (for large enough $M$)
\[
\mbox{$r(\uj) \approx r(j_{\pi(1)}\dots j_{\pi(k)})
\,\,\, \forall \uj\in [K]^k$ and $\pi\in \mS_{[k]}$,}
\]
or, since adjacent transpositions generate $\mS_{[k]}$
(here ``adjacent'' is just notationally convenient),
\beq{rjrj}
\mbox{$r(\uj) \approx r(j_{1}\dots j_{s-1},j_{s+1},j_s,j_{s+2}\dots j_k)
\,\,\, \forall \uj\in [K]^k$ and $s\in [k-1]$}.
\enq

\nin
\emph{Proof of \eqref{rjrj}.}
We may assume $M=km$ and let $I_1\cup\cdots\cup I_k$  partition $[M]$ 
into intervals of size $m$.  Set
\[
a(\ell,j) = m^{-1} \sum_{i\in I_\ell} \mbone_{\{X_i=j\}}
\]
and notice that \eqref{r(j)} implies 
\beq{Eprod}
\E \prod_{\ell=1}^k a(\ell,j_\ell) \approx r(\uj)  \,\,\, \forall\uj\in [K]^k.
\enq
Fix $\uj$ and set $\psi =\prod\{ a(\ell,j_\ell):\ell\in [k]\sm\{ s,s+1\}\}$.
A second application of \eqref{r(j)} gives, for any $j$,
\beq{Esquare}
\E [\psi \cdot(a(s,j)-a(s+1,j))^2] =O(\eta).
\enq
[\emph{Because (briefly)}:  The l.h.s.\ is
\[
\mbox{$m^{-k}\left\{\sum'\E [\psi \cdot \mbone_{\{X_i=X_{i'}=j\}}]-
\sum''\E [\psi \cdot \mbone_{\{X_i=X_{i'}=j\}}]\right\},$}
\]
where the $(i,i')$'s in $\sum'$ range over $I_s^2\cup I_{s+1}^2$ and those in
$\sum''$ over $(I_s\times I_{s+1})\cup (I_{s+1}\times I_s)$.
Of the $m^k$ terms in each sum, all but the diagonal terms in $\sum'$ are 
within $O(\eta)$ of $r(j_{1}\dots j_{s-1},j,j,j_{s+2}\dots j_k)$.
This gives \eqref{Esquare} if we pretend \emph{all} terms are non-diagonal; 
and if, say, $m>1/\eta$, the effect of the $2m^{k-1}$
diagonal terms is negligible.]

\mn

For \eqref{rjrj} we now have, 
setting $\uj'=(j_{1}\dots j_{s-1},j_{s+1},j_s,j_{s+2}\dots j_k)$,
writing $\|f\|$ for $\sqrt{\E f^2}$, and with the initial ``$\approx$'' given by \eqref{Eprod}:
\begin{eqnarray*}
|r(\uj)-r(\uj')|&\approx&
|\E\left[\psi\cdot \left\{a(s,j_s)a(s+1,j_{s+1}) -
a(s,j_{s+1})a(s+1,j_s)\right\}\right]|\\
&\leq &
\E| \psi\cdot \left\{a(s,j_s)(a(s+1,j_{s+1}) -a(s,j_{s+1}))\right\}|\\
&& \,\,\,\,\,\,\,\,\,\,\,\,\,\,\,\,\,\,\,\,\,\,\,\,\,\,\,\,\,\,\,\,\,
+ \E| \psi\cdot \left\{a(s,j_{s+1})(a(s,j_s) - a(s+1,j_s))\right\}|\\
&\leq &
\|\sqrt{\psi}\cdot   a(s,j_s)\| \|\sqrt{\psi} (a(s+1,j_{s+1}) -a(s,j_{s+1})) \|\\
&& \,\,\,\,\,\,\,\,\,\,\,\,\,\,\,\,\,\,\,\,\,\,\,\,\,\,\,\,\,\,\,\,\,
+ 
\|\sqrt{\psi}\cdot   a(s,j_{s+1})\| \|\sqrt{\psi} (a(s,j_s) - a(s+1,j_s)) \|\\
&\leq &
\|\sqrt{\psi} (a(s+1,j_{s+1}) -a(s,j_{s+1})) \|
+ 
\|\sqrt{\psi} (a(s,j_s) - a(s+1,j_s)) \|\\
&=&O(\sqrt{\eta})
\end{eqnarray*}
(where the last ``$\leq$'' holds because $\psi$ and $a(\cdot,\cdot)$ are 
in $[0,1]$), and the final bound follows from \eqref{Esquare}.

\qed

%\newpage
\section{Examples}\label{Examples}

Here we fill in examples promised in remarks following Theorem~\ref{T1/e}
and Lemma~\ref{LWin}, 
the first showing that Theorem~\ref{T1/e}(a) is locally almost optimal, 
and the second that even a large collection of elements of $P$ of equal average heights
need not contain a pair with balance much better than the $1/e$ of \eqref{1/e}.

We need a pair of (fairly) standard definitions:
the \emph{parallel} and \emph{series sums} of posets $P$ and $Q$, 
denoted $P+Q$ and $P \oplus Q$,  
are, respectively, the disjoint union of $P$ and $Q$ (with no added relations), and the 
poset obtained from $P+Q$ by setting $P<Q$.
We use $C_k$ for the chain of size $k$.

\begin{Ex}\label{MaxT4.1}
For any $\eps>0$, there is a poset $P$ with antichain $I$ satisfying $|I|\gg n^{1/2-\eps}$ and 
\beq{gdxy}
\mbox{$\delta_{xy}\le \eps$ for all distinct $x,y\in I$.}
\enq
\end{Ex}
\nin
\emph{Proof.} For $I$ an antichain of $P$, call $(P,I)$ \emph{good} if \eqref{gdxy} holds. 
It suffices to show that if $(P,I)$ is good and $m:=|P|$ is sufficiently large, then there is a good $(Q,J)$ 
with $|J|=2|I|$ and $|Q|= (4+2\eps)m$. 
To see this, let $(P',I')$ and $(P'',I'')$ be copies of 
$(P,I)$, and 
\[
(Q,J)=([P'\oplus A'] +[A''\oplus P''], I'\cup I''),
\]
where $A'$ and $A''$ are antichains of size $(1+\eps)m$, 
and notice that $(Q,J)$ is good (for large enough $m$), 
since $P'\prec P''$ with probability $1-2^{-\gO(m)}$.

\qed

We note that a similar construction gives $P$ and $I$ with $|I|\gg 1$ and $\delta_x\le \eps$
$\forall x\in I$.  (It had at one time seemed plausible that a large antichain must 
contain an $x$ with $\gd_x > 1/2-o(1)$.)

\begin{Ex}\label{MaxT2.3}
For any $\eps>0$ and $k$ there are a poset $P$ and $x_1\dots x_k\in P$
with $H(x_i)=1/2$ $\forall i$ and $\delta_{x_ix_j}\le 1/e+\eps$ $\forall i\neq j$.
\end{Ex}
\nin
The justification here is longer than for Example~\ref{MaxT4.1} and will skip a few routine verifications.
We first note a pair of easy observations, whose proofs we omit.
Given $P$,
%of size $n$, 
we now use $\tdf$ for $f/(|P|+1)$.
\begin{prop}\label{prop2}
{\rm (a)}
For any $\eta,\gd>0$ there is an $L$ such that if $|P|\ge L$ and
$x\in P$, then
\[
\pr(|F(x)-H(x)|>\gd)<\eta .
\]
{\rm (b)}  
For any $\eta>0$ there is a $C$ such that for any $P$ and
$x\in P$, 
$\pr(|F(x)-\tdf(x)|>C/\sqrt{|P|})<\eta .$

%For any $\eta>0$ there is a $C$ such that for any $P$ of size $n$ and $x\in P$, 
%$\pr(|F(x)-\tdf(x)|>C/\sqrt{n})<\eta .$

\end{prop}

%\emph{Proof.}
%This will follow easily from the next two observations.
The main point for Example~\ref{MaxT2.3} is the following construction.

\begin{prop}\label{prop1}
For any $\eta>0$ and $L$ there are $\gd>0$, $P$ and $x\in P$ with $|P|> L$, $H(x)=1/2$ and
\[
\pr(F(x)> 1/2-\gd)<1/e +\eta.
\]

\end{prop}

Before proving this 
%Proposition~\ref{prop1} 
we show that it gives
Example~\ref{MaxT2.3}. 
Let $(P_1,x_1,\gd_1)$ be as in the proposition
with $\eta=\eps/2$;
for $i=2\dots k$, let $(P_i,x_i,\gd_i)$ be as in the proposition with $\eta =\eps/e2$
and $|P|$ large enough
that $\pr(F(x_i)< 1/2- \gd_j)< \eps/2$ $\,\,\forall j<i$ (see Proposition~\ref{prop2}(a));
and let $P=P_1+\cdots+P_k$. Then for $1\le j< i\leq k$,
\[
\pr(F(x_j)>F(x_i))\le \pr(F(x_j)> 1/2-\gd_j)+\pr(F(x_i)< 1/2-\gd_j)< 
(1/e+\eps/2)+\eps/2= 1/e+\eps.
\]\qed

\nin
\emph{Proof of Proposition~\ref{prop1}.}
Since this discussion will involve several posets we now use $F_P$ rather than simply $F$,
and similarly for other quantities.
Requirements for $r,a,k,l$ and $\gd$ will appear below.

Let $r\gg 1$ and
let $Q=C_r$, with minimum element $x$.
Then $H_Q(x) = 1/(r+1)$ (see \eqref{H}),
\beq{prFQ}
\pr(F_Q(x)>\ga) = (1-\ga)^r \,\,\,\forall \ga\in [0,1].
\enq
In particular $\pr(F_Q(x)> H_Q(x))\approx 1/e$, and
the goal of the rest of the construction is to preserve this behavior (more or less)
while shifting $H(x)$ to 1/2.

Let
$R= Q+C_a$ and $S=C_k\oplus R \oplus C_l$ ($S$ will be the $P$ of the proposition), 
and set $m = |S| = k + r+a+l$.  Then
\[
H_R(x) =H_Q(x) = 1/(r+1),
\]
\[
h_S(x) = k+ h_R(x) = k +(r+a+1)/(r+1) = k + (r+a+1)H_Q(x).
\]
Choose (with elaboration below) $r\gg 1$, $a\gg r^2$ 
and $k,l < 3a$, satisfying
\[
H_S(x) = h_S(x)/(m+1) = [k +(r+a+1)/(r+1)]/(m+1) =1/2
\]
(this is easily seen to be possible), noting that then
\beq{m+1}
[(m+1)/2-k ]/(a+r+1) = 1/(r+1).
\enq

Now aiming for the inequality in the proposition, we have 
\beq{prFSx}
\pr(F_S(x) > 1/2 -\gd) < \pr (\tdf_S(x) > 1/2-2\gd) +\pr(|F_S(x)-\tdf_S(x)|>\gd)
\enq
and 
\begin{eqnarray}\label{prf's}
\pr (\tdf_S(x) > 1/2-2\gd) &=& \pr(f_S(x) > (m+1)(1/2-2\gd))\nonumber\\
&=& \pr(f_R(x) > (m+1)/2-k -2(m+1)\gd )
\nonumber\\
&\le & 
%\pr \left(F_R(x) > \tfrac{(m+1)/2-k -2(m+1)\gd}{a+r+1} -\gd\right)
\pr \left(F_R(x) > \tfrac{1}{r+1} -\left[\tfrac{2(m+1)}{a+r+1}+1\right]\gd\right)
+\pr(|\tdf_R(x)-F_R(x)|>\gd)
\end{eqnarray}
(where \eqref{prf's} uses \eqref{m+1}).
Now choosing $\gd$ satisfying $1/\sqrt{a} \ll \gd \ll 1/r$, we find that 
the last terms in
\eqref{prFSx} and \eqref{prf's} are small (say less than $\eta/3$) by Proposition~\ref{prop2}(b)
(so we want $\gd> C/\sqrt{a}$ for a suitable $C$ depending on $\eta$),
while, by \eqref{prFQ}, taking $r$ at least about $1/\eta$ and
$\gd< c/r$ for a suitably small $c>0$ makes the first term in
\eqref{prf's} less than $e^{-1}+\eta/3$.

\qed

\section{Problems}\label{Conjectures}

We have already seen some of these (including, in addition to 
the venerable conjectures of Section~\ref{Intro},
the conjectures and questions of Section~\ref{Parameters}, and Conjecture~\ref{CWinVar}).  
Here we collect a few more of the many questions that arose in the course of this work.
Though mostly originating with balance problems,
these seem to us to be of independent interest, and serve as a 
(maybe unnecessary) reminder of how poorly we understand this whole subject.
(We haven't thought equally hard about all of these questions, and it could certainly be that
some are either wrong or relatively easy.)

\nin

We begin with a possibility that strengthens Theorem~\ref{Tgapw} but is weaker than
Conjecture~\ref{A1}.
\begin{conj}\label{CgapwP}
For any P, $\gap(P) = O(w(P))$.
\end{conj}
\nin
This
%, with Theorem~\ref{TKS}, 
would have a nice corollary which, 
amazingly, we don't know to be true:
\beq{isox}
\mbox{\emph{if there is an $x$ in $P$ 
%contains some $x$ 
comparable to no other element, then 
$P$ satisfies \eqref{1/2}.}}
\enq
(It is \emph{not} true that such an $x$ must itself belong to a pair witnessing \eqref{1/2}.)

\nin
\emph{Derivation of \eqref{isox} from Conjecture~\ref{CgapwP}} (briefly):
If \eqref{1/2} fails for $P$, then Theorem~\ref{TKS}(c) gives $w(P) < \gz n$ with $\gz=o(1)$.
But then, by Conjecture~\ref{CgapwP}, $\gap(P) =O(\gz n)$, so there is a $y\in P\sm x$
with $h(y) = (1\pm O(\gz))n/2 = (1\pm o(1))n/2$; and $\gd_{xy} > 1/2-o(1)$ for any such $y$.
\qed

\nin
%\textbf{C.}

For a chain
$C=\{y_0<\cdots <y_m \}$, define the \emph{gap} of $C$ to be 
\beq{gapC}
\gap(C) = \max\{h(y_1),h(y_2)-h(y_1)\dots h(y_m)-h(y_{m-1}), n+1-h(y_m)\},
\enq
and for $x\not \in C$ (and $\eps \in [0,1]$), say $(x,C)$ 
is $\eps$-\emph{diffuse} if 
\beq{diffuse}
\max\{p(x\prec y_1),p(y_1\prec x\prec y_2)\dots p(y_{m-1}\prec x\prec y_m),p(x\succ y_m)\} < \eps.
\enq

Though we're not ready to conjecture, 
as far as we know even the following very strong versions of Conjecture~\ref{CgapwP}
are possible.
\begin{question}\label{chain.gaps}
Is it true that any P contains a chain C with $\gap(C)=O(w(P))$?
Could it be that any longest chain has this property?
\end{question}

Notice that if \eqref{diffuse} holds, then $\gd_{xy}> 1/2-\eps$ for some $y\in C$.
(This simple observation underlies Theorem~\ref{TKS}(d) 
[see \eqref{xinsome}], as well as, e.g., Linial's proof of \eqref{1/3} for $P$ of
width 2 \cite{Linial}.)
Thus, the following speculations would substantially strengthen Theorem~\ref{Twsig}.

\begin{question}\label{CTwin'}
If $w(P)$ is fixed and $\gs(x)\gg 1$, must
$(x,C)$ be $o(1)$-diffuse for some chain $C$?
Could this be true whenever $C$ is a longest chain of $\Pi(x)$?
\end{question}
Before moving on from gaps, we mention the following explicit proposal for
Conjecture~\ref{Qparams}, saying that an $X$ witnessing large $\tau(P)$ can
be found among the elements ``bordering''
any large gap.
\begin{conj}\label{Cgaptau}
If $D$ is an ideal of $P$, $U=P\sm D$, 
%$A=\max (D)$ and $B=\min (U)$.
and $\min_{x\in U}h(x) - \max_{x\in D}h(x) \gg 1$ then there is an $X$ 
contained in either $\max (D)$ or $\min (U)$ with
$|X|^{-1}H(\gss|_X)\gg 1$.
\end{conj}

%\nin
%\textbf{D.}

\nin
%\emph{Perturbations}

In thinking about the present material,
it has often seemed of interest to understand effects of small changes in $P$, 
e.g.\ as in the next pair of suggestions (in which $h_P$ etc.\ have the natural meanings).

\begin{conj}\label{increase.h}
For any P and x:

\nin
{\rm (a)}  
$\,\,\,\sum_{y\neq x}\max\{h_P(y)-h_{P-x}(y),0\}\leq n-1.$

\nin
{\rm (b)}  
$\,\,\,\sum_{y\neq x} |H_P(y)-H_{P-x}(y)|\leq (n-1)/(n+1).$
\end{conj}

\nin
As far as we know, (a) could still hold with the summand replaced by $|h_P(y)-h_{P-x}(y)|$.  
%We also believe $H_P(y)-H_{P-x}(y)\ge 0$ if $x<y$, but only know this if $x\in \min (P)$,
%where it doesn't require $x<y$.

Both parts of Conjecture~\ref{increase.h} 
are exact if $x=\hO$ (or $\hone$), and 
true if $x$ is minimal (or maximal).
Here it's important to consider \emph{global} effects, since
individual summands can be of the same order as the conjectured bounds
(e.g.\ when $x<y$ is the only relation involving $x$ or $y$).
On the other hand, we do expect smaller local effects in 
many natural situations; for example:
\begin{conj}\label{Hdeletex}
We have $|h_P(y)-h_{P-x}(y)|=O(1)$ assuming any one of:

\nin
{\rm (a)}  there is a chain $x=x_0<x_1<\cdots < x_k=y$ with $k=\gO(n)$;

\nin
{\rm (b)}  the distance between $x$ and $y$ in the cover graph of $P$ is $\gO(n)$;

\nin
{\rm (c)}  $\pi(x) =O(1)$.

\end{conj}

\nin
(The \emph{cover graph} is the graph on $P$ with $x$ and $y$ adjacent 
if one of them covers the other.)
As far as we know, the weaker $\win(x)=O(1)$ could also be enough in (c).

\nin

Our last two questions are more miscellaneous.

As said following Theorem~\ref{TKS'}, $w(P)>\log n$ doesn't
imply \eqref{1/k!}
(i.e.\ $\gd_k(P)\ra 1/k!$) even for $k=3$.

\begin{question}\label{Qgdk}
For how large a $w=w(n)$ can one find $n$-element, width $w$ posets violating \eqref{1/k!}
for $k=3$?
\end{question}
\nin
(In the other direction we know only Theorem~\ref{TKS'}(c).
At this writing we don't know that width $\go(\log n)$
%$w(P)\gg \log n$ 
doesn't imply $\gd_3\ra 1/6$, though we \emph{can}
give examples with width $n^{1/2-o(1)}$ that violate \eqref{1/k!} for $k=4$.)
%(which suggests the possibility of some hierarchy).}

\mn

Finally, we return to the frustrating Conjecture~\ref{CWinVar}
and offer the following precise version, which is 
tight whenever $\Pi(x)\not\sim P\sm\Pi(x)$, and is, 
at least, correct for series-parallel posets (defined e.g.\ in \cite{Lawler}).

\begin{conj}\label{CWinVar'}
For any P and x, 
%with $\Win(x) = W$ and $H(x)=H$,
\beq{VarF'}
\Var (F(x))\leq \frac{\Win(x)  H(x)(1-H(x))}{2+\Win(x)}.
\enq
\end{conj}
\nin
In particular, when $x\in \min(P)$ (so $\Win(x) = H(x)/2$), \eqref{VarF'} becomes 
the appealing (also unproved)
\[
\Var (F(x))\leq \frac{H^2(x)(1-H(x))}{1+H(x)}.
\]

\end{document}